\documentclass[final,onefignum,onetabnum]{siamonline250211}

\usepackage{amsmath,amsfonts}
\usepackage{graphicx}
\usepackage{graphbox}
\usepackage{subfigure}
\usepackage{epstopdf}
\usepackage{algpseudocode}
\ifpdf
  \DeclareGraphicsExtensions{.eps,.pdf,.png,.jpg}
\else
  \DeclareGraphicsExtensions{.eps}
\fi

\usepackage{enumitem}
\setlist[enumerate]{leftmargin=.5in}
\setlist[itemize]{leftmargin=.5in}

\newcommand{\ud}{\mathrm{d}}
\newcommand{\e}{\mathrm{e}}
\renewcommand{\i}{\mathrm{i}}

\newsiamremark{remark}{Remark}
\newsiamremark{hypothesis}{Hypothesis}
\crefname{hypothesis}{Hypothesis}{Hypotheses}
\newsiamthm{claim}{Claim}

\headers{Estimation of a Gas Diffusion Coefficient}{I. Viviano}

\title{Estimation of a Gas Diffusion Coefficient by Fitting Molecular Dynamics Trajectories to Finite-Difference Simulations}

\author{Isaac Viviano\thanks{(\email{iviviano@oberlin.edu})}}

\dedication{\small\textit{Project advisor: Robert Krasny\thanks{University of Michigan, Ann Arbor, MI}}}

\usepackage{amsopn}

\ifpdf
\hypersetup{
  pdftitle={Estimation of a Gas Diffusion Coefficient from Molecular Dynamics Simulations},
  pdfauthor={I. Viviano}
}
\fi

\begin{document}

\maketitle

\begin{abstract} 
A procedure is presented to estimate the diffusion coefficient
of a uniform patch of argon gas in a 
uniform background of helium gas.
Molecular Dynamics (MD) simulations 
of the two gases interacting through the 
Lennard-Jones potential 
are carried out using the LAMMPS software package.
In addition,
finite-difference (FD) calculations are used to solve the 
continuum diffusion equation for the
argon concentration with a given diffusion coefficient.
To contain the computational cost
and facilitate data visualization,
both MD and FD computations were done in two space dimensions.
The MD argon trajectories were binned to the FD grid,
and
the optimal diffusion coefficient
was estimated by minimizing the 
difference between the binned MD data 
and the FD solution with a
nonlinear least squares procedure
(Levenberg–Marquardt algorithm).
Numerical results show
the effect of the MD binning parameter
and
FD grid spacing.
The estimated diffusion coefficient is compared to an experimental measurement.
\end{abstract}


\section{Introduction} 

The diffusion of a solute gas 
(for example argon) through a solvent gas
(for example helium)
can be quantified by a diffusion coefficient~$D$.
The coefficient can be measured
experimentally by gas chromatography, which separates a mixture of gases into pure components as it travels through a column so that $D$
can be calculated from the peak widths 
of the mass spectrum~\cite{wasik_1969_measurements}.
The diffusion coefficient can also be calculated theoretically
from the expression
\begin{equation} \label{eq:coeff_formula}
    D=\frac{1}{6}\lim_{t\to\infty}\left\langle \frac{\ud}{\ud t}(\mathbf{r}(t)-\mathbf{r}(0))^2\right\rangle,
\end{equation}
where $\mathbf{r}(t)$ is the trajectory of a gas particle 
and 
$\langle \cdot\rangle$ represents the average over all diffusing particles~\cite{melzer}.
The derivative of the square displacement describes how fast each particle spreads out,
so for sufficiently large simulation times, 
the ensemble average in \cref{eq:coeff_formula} quantifies the overall diffusion rate. 
One way to generate the particle trajectories for \cref{eq:coeff_formula} is with 
Molecular Dynamics (MD) simulations,
which provide the trajectories of a collection of particles evolving by Hamiltonian dynamics
in a classical interparticle force field.
Here, we present an alternative theoretical approach to \cref{eq:coeff_formula} for determining the 
gas diffusion coefficient.

Fick's laws give a continuum description
of diffusion,
resulting in a partial differential equation
for the time-evolution of the solute concentration~\cite{Fick}.
Our approach estimates the the diffusion coefficient~$D$
by minimizing the discrepancy
between numerical finite-difference (FD) solutions of the 
diffusion equation 
and MD simulations of the solute concentration.
We anticipate that this procedure 
can be applied to estimate coefficients
in a differential equation for other physical systems 
described by a continuum model such as phase separation of polymers.

Our MD simulations are carried out using LAMMPS
(Large-scale Atomic/Molecular Massively Parallel Simulator),
an open-source classical MD package focused on materials modeling,
which was developed by researchers at
Sandia National Laboratories and collaborators~\cite{LAMMPS}.
We fix values of the state variables (pressure and temperature)
and conduct equilibrium MD simulations
where the gas particles follow Hamiltonian dynamics of 
pairwise Lennard-Jones interactions.
Since MD simulations can incorporate the effect 
of varying the state variables,
our procedure could be applied to 
estimate the dependence of the diffusion coefficient
on these variables. 

In general,
MD and FD calculations may be applied to 
one-, two-, and three-dimensional systems,
and
throughout this work the spatial dimension is denoted $d$.
Similarly, 
LAMMPS provides the functionality to conduct MD simulations in one-, two-, or three-dimensions.
While the models are equally valid for 
any dimension,
the estimation procedure is implemented for $d=2$ 
to contain the
computational cost and facilitate data visualization.

The rest of the article is organized as follows.
The derivation of the diffusion equation 
is given in~\cref{sec:HE},
and
the FD methods are 
presented in~\cref{ssec:heat_fd}.
The MD approach is described in~\cref{sec:MD}, 
and the MD and FD simulation results are compared in~\cref{sec:res}.
The diffusion coefficient estimation algorithm is presented in~\cref{sec:estim},
and the estimates are given in~\cref{sec:est_res}.
The work is summarized in~\cref{sec:sum}.

%%%%%%%%%%%%%%%%%%%%%%%%%%

\section{Diffusion equation} \label{sec:HE}

This section introduces the continuum model for 
gas diffusion, 
explaining its physical basis
and its solution by Fourier series.
The gas concentration is a scalar function
of space and time,
$u(\mathbf{x},t),~\mathbf{x} \in \mathbb{R}^d$.
For a system of gas particles subject to a 
non-zero concentration gradient, 
random motion tends to reduce the 
gradient over time.
The net movement of particles from regions of higher concentration to lower concentration is called diffusion
and
these concepts are quantified by Fick's laws.
While our numerical methods do not make use of the analytical solution,
it prescribes important physical properties of the differential equation that the finite-difference methods should replicate. 
In particular,
the periodic boundary conditions ensure that mass is conserved and energy decreases in time. 

%%%%%%%%%%%%%%%%%%%%%%%%%%

\subsection{Derivation from Fick's laws}

Fick's first law says that the 
diffusion flux ${\bf J}({\bf x},t)$
is proportional to the 
concentration gradient $\nabla u({\bf x},t)$
with proportionality constant 
given by a scalar diffusion coefficient $D$,
\begin{equation} \label{eq:fick1}
\mathbf{J}({\bf x},t) = -D\nabla u({\bf x},t).
\end{equation}
In general
the diffusion coefficient depends on
many factors such as the gas temperature and pressure,
but in this work it is assumed to be a given constant.
Conservation of mass states that the rate of change
of the mass contained in a volume $V$ 
is determined by the mass flux through the 
surface $S =\partial V$ bounding the volume~$V$,
\begin{equation}
\frac{d}{dt}\int_V u({\bf x},t) \ud\mathbf{x} +
\oint_{S}\mathbf{J}(\mathbf{x},t)\cdot 
\ud {\bf S} =0.
\end{equation}
Application of the divergence theorem 
converts the surface integral into a volume integral
which yields the differential form of 
mass conservation,
\begin{equation} \label{eq:differential_mass}
\frac{\partial}{\partial t}u(\mathbf{x},t)+\nabla \cdot\mathbf{J}(\mathbf{x},t)=0.
\end{equation}
Then substituting \cref{eq:fick1} into \cref{eq:differential_mass}
yields Fick's second law,
\begin{equation} \label{eq:HE}
\frac{\partial}{\partial t }u(\mathbf{x},t) = 
D\nabla^2u(\mathbf{x},t), 
\end{equation}
a partial differential equation
for the time-evolution of the 
gas concentration
known as the diffusion equation~\cite{Fick}.
The operator $\nabla^2 = \nabla \cdot \nabla$ 
on the 
right side of \cref{eq:HE} is the Laplacian.
This work considers the diffusion equation
with periodic boundary conditions 
on the $d$-dimensional unit cube $\Omega=[0,1]^d$
and
a given initial condition,
\begin{subequations} \label{eq:bc_ic}
\begin{align}
&u(\mathbf{x},t)\big|_{x_i=0} = 
u(\mathbf{x},t)\big|_{x_i=1}, \quad
\frac{\partial}{\partial x_i} u(\mathbf{x},t)\big|_{x_i=0} = 
\frac{\partial}{\partial x_i} u(\mathbf{x},t)\big|_{x_i=1}, \quad
i=1,\ldots, d, 
\label{subeq:bc} \\
&u(\mathbf{x},0) = u_0(\mathbf{x}).
\label{subeq:ic}
\end{align}
\end{subequations}
Note that \cref{eq:HE} 
is also known as the heat equation since
it models heat flow when
$u$ represents the material temperature and $D$ is its thermal conductivity. 
The next subsection derives an analytic expression
for the solution of the diffusion equation.

%%%%%%%%%%%%%%%%%%%%%%%%%%%%%%%%%%

\subsection{Analytical solution by Fourier Series} \label{ssec:analytical}

For the $d$-dimensional unit cube 
with periodic boundary conditions, 
the solution of the diffusion equation \cref{eq:HE} can be expressed as a Fourier series~\cite{pdebook},
\begin{equation} \label{eq:heat_solution}
    u(\mathbf{x}, t)=\sum_{\mathbf{m}\in\mathbb{Z}^d}\widehat{u}_\mathbf{m}(t)\e^{2\pi\i\mathbf{m}\cdot\mathbf{x}},\quad 
    \widehat{u}_\mathbf{m}(t)=\int_{\Omega}u(\mathbf{x},t)~\e^{-2\pi\i \mathbf{m\cdot x}}~\ud \mathbf{x},
\end{equation}
where $\widehat{u}_\mathbf{m}(t)$ 
denotes the ${\bf m}$-th 
Fourier coefficient of $u({\bf x},t)$ 
with wavenumbers over the $d$-dimensional integer lattice ${\bf m} \in \mathbb{Z}^d$.
From substituting \cref{eq:heat_solution} into \cref{eq:HE} we obtain the initial value problem,
\begin{equation}
\label{eq:fourier_coefficient_ODE}
\frac{\ud}{\ud t} \widehat u_{\bf m}(t) =
-4\pi^2|\mathbf{m}|^2D \widehat u_{\bf m}(t),
\quad\widehat{u}_\mathbf{m}(0) = 
\int_\Omega u_0(\mathbf{x})~\e^{-2\pi\i\mathbf{m\cdot x}}~\ud \mathbf{x},
\end{equation}
for the Fourier coefficients.
Its solution yields the time dependence,
\begin{equation} \label{eq:time}
\widehat{u}_\mathbf{m}(t) =
\widehat{u}_\mathbf{m}(0)
\e^{-4\pi^2|\mathbf{m}|^2Dt},
\end{equation}
where $\widehat{u}_\mathbf{m}(0)$ is the 
$\mathbf{m}$-th Fourier coefficient of the initial condition $u_0$.
Several properties of the solution 
can now be derived.

First,
setting ${\bf m} = {\bf 0}$ in~\cref{eq:time}
yields 
$\widehat{u}_{\bf 0}(t) = 
\widehat{u}_\mathbf{0}(0)$,
which implies that the total mass of the system
is conserved in time,
\begin{equation} \label{eq:mass_con}
\int_\Omega u({\bf x},t)~\ud\mathbf{x} = 
\int_\Omega u_0({\bf x})~\ud\mathbf{x}.
\end{equation}
In addition,
non-constant modes ${\bf m} \ne {\bf 0}$
decay exponentially in time, 
and hence the concentration $u({\bf x},t)$ 
converges to its average value
$\int_\Omega u_0(\mathbf{x})~\ud\mathbf{x}$
as $t\to\infty$. 

Another important property of solutions
to the diffusion equation with
periodic boundary conditions
is that their energy is non-increasing in time.
The energy is given by the 
square of the 2-norm of the solution,
\begin{equation}
E(t) = \|u(\cdot,t)\|_2^2 =
\int_\Omega|u(\mathbf{x},t)|^2\ud\mathbf{x} =
\sum_{\mathbf{m}\in\mathbb{Z}^d}|
\widehat{u}_\mathbf{m}(t)|^2,
\end{equation}
where the last equality follows from
Parseval's relation which states that the 
2-norm of a periodic function equals
the discrete 2-norm of its Fourier coefficients.
Using the ODE for the Fourier coefficients~\cref{eq:fourier_coefficient_ODE},
it follows that
\begin{equation}
E^\prime(t) =
2\sum_{\mathbf{m}\in\mathbb{Z}^d}
\widehat{u}_\mathbf{m}(t)
\frac{\ud}{\ud t} \widehat u_{\bf m}(t) =
-8\pi^2D\sum_{\mathbf{m}\in\mathbb{Z}^d}
|\mathbf{m}|^2|\widehat u_{\bf m}(t)|^2 \le 0,
\end{equation}
which implies that the energy is
non-increasing in time.
The next section discusses 
finite-difference methods which are used
to solve the diffusion equation numerically.

%%%%%%%%%%%%%%%%%%%%%%%%%%%%%%%%%%%%%%%%%%%%%%%

\section{Finite-Difference method} \label{ssec:heat_fd}

In a finite-difference method
the temporal and spatial domains are discretized
and
the derivatives in the equation
are replaced
by finite-difference approximations
that enable the numerical solution to be found
by stepping forward in time.
This work considers two such methods
for solving the diffusion equation~\cref{eq:HE}
with 
periodic boundary conditions and initial condition~\cref{eq:bc_ic}.
A time step $k$ 
and space step $h=1/N$ are chosen,
where $N$ is a positive integer.
The associated spatial grid is
$\Omega_h = h[N]^d$,
where $[N] = \{0,\ldots,N-1\}$,
and
the solution $u({\bf x},t)$
is approximated at grid points
${\bf x}_{\bf j} = h{\bf j}$
and time steps $t_n = nk$,
\begin{equation} \label{eq:seq_approx}
u_\mathbf{j}^n \approx u({\bf x}_{\bf j},t_n),
\quad {\bf j} = (j_1,\ldots,j_d) \in [N]^d,
\quad n=0,\ldots,n_{\rm max}.
\end{equation}
Next, we consider the spatial discretization
of the Laplacian.

%%%%%%%%%%%%%%%%%%%%%%%%%%

\subsection{Central difference approximation in space}

This work employs a 
central difference approximation
for the spatial discretization of the
Laplace operator,
yielding the discrete Laplacian, 
\begin{equation}
\nabla^2_hu_\mathbf{j}^n=\frac{1}{h^2}\sum_{i=1}^d (u_{\mathbf{j}+\mathbf{e}_i}^n-2u_\mathbf{j}^n+u_{\mathbf{j}-\mathbf{e}_i}^n),
\quad {\bf j} \in [N]^d,
\end{equation}
where ${\bf e}_i$ is 
the $i$-th standard basis vector in $\mathbb{R}^d$.
The vector $u_{\mathbf j}^{n}$ is extended
to an $N$-periodic signal,
$u_\mathbf{j}^n=u_{\mathbf{j}+N\mathbf{e}_i}^n$
for ${\bf j}\in{\mathbb Z}^d, i=1,\ldots,d$.
This justifies the discrete Fourier series expansion,
\begin{equation} \label{eq:as_DFT}    u_{\mathbf{j}}^{n} =
\frac{1}{N^d}\sum_{\mathbf{m}\in[N]^{d}}
\widehat{u}_{\mathbf{m}}^{n}\e^{2\pi\i
\mathbf{m \cdot j}/N}, \quad
\widehat{u}_{\mathbf{m}}^{n} =
\sum_{\mathbf{j}\in[N]^{d}}
u_{\mathbf{j}}^{n}\e^{-2\pi\i 
\mathbf{m \cdot j}/N},
\end{equation} 
where 
$\widehat{u}_\mathbf{m}^n$ is the 
${\bf m}$-th discrete Fourier coefficient 
of the numerical solution $u_{\bf j}^{n}$.
Using properties of the exponential function, 
we can understand how the discrete Laplacian acts in Fourier space,
\begin{equation} \label{eq:mult}
\begin{split}
\nabla^2_{h}u_{\mathbf{j}}^{n}&=\frac{1}{h^2}\sum_{i=1}^{d} \frac{1}{N^d}\sum_{\mathbf{m}\in[N]^{d}}\widehat{u}_{\mathbf{m}}^{n}(\e^{2\pi \i \mathbf{m}\cdot(\mathbf{j} + \mathbf{e}_{i})/N}-2\e^{2\pi\i\mathbf{m\cdot j}/N}+\e^{2\pi\i\mathbf{m}\cdot(\mathbf{j}-\mathbf{e}_{i})/N})\\		&=\frac{1}{N^d}\sum_{\mathbf{m}\in[N]^{d}}\sum_{i=1}^{d}\widehat{u}_\mathbf{m}^n\frac{1}{h^2}( \e^{2\pi\i m_i/N}+\e^{-2\pi\i m_i/N}-2 )\e^{2\pi\i\mathbf{m\cdot j}/N}\\
		&= \frac{1}{N^d}\sum_{\mathbf{m}\in[N]^{d}}\widehat{u}_{\mathbf{m}}^{n}\frac{1}{h^2}\sum_{i=1}^{d}(\e^{2\pi\i m_i/N}+\e^{-2\pi\i m_i/N}-2)\e^{2\pi\i \mathbf{m\cdot j}/N} \\	&=\sum_{\mathbf{m}\in[N]^{d}}\lambda_\mathbf{m}\widehat{u}_\mathbf{m}^n\e^{2\pi\i\mathbf{m\cdot j}/N}, \quad
        \lambda_{\mathbf{m}} 
        = \frac{1}{N^d}\frac{1}{h^2}\sum_{i=1}^d(\e^{2\pi\i m_i/N}+\e^{-2\pi\i m_i/N}-2),
  \end{split}
\end{equation}
where $\lambda_\mathbf{m}$ are the eigenvalues 
and
$(\e^{2\pi\i \mathbf{m}\cdot\mathbf{j}/N} : \mathbf{j} \in [N]^d)$
are the corresponding
eigenvectors
of the 
discrete Laplacian with 
periodic boundary conditions.
\Cref{eq:mult} illustrates the principle that
applying the spatial operator in physical space is equivalent to multiplication in Fourier space.
An alternative expression for the eigenvalues is given by
\begin{equation} 
\label{eq:eigenvalue}
    \lambda_{\mathbf{m}} 
    = \frac{2}{h^{2}}\sum_{i=1}^{d}
    (\cos(2\pi m_{i}/N)-1)
    = - \frac{4}{h^{2}}\sum_{i=1}^{d}\sin^{2}(\pi m_{i}/N), \quad {\bf m} \in [N]^d,
\end{equation}
and since $0\le\sin^{2} (\pi m_i/N)\le 1$,
the eigenvalues are nonpositive and bounded,
\begin{equation} \label{eq:apriori}
-4d/h^2 \le \lambda_\mathbf{m} \le 0.
\end{equation}

Next we consider the stability of the difference scheme.
Recall that the energy of the solution of the diffusion equation~\cref{eq:HE}
does not increase in time,
and
the difference scheme should have the same property.
This is expressed as the discrete stability condition,
\begin{equation}
\label{eq:discrete_stability}
\|{\bf u}^{n+1}\|_2 \le
\|{\bf u}^n\|_2, \quad
\|{\bf u}^n\|_2^2 =
\sum_{\mathbf{j}\in[N]^d }|u_\mathbf{j }^{n}|^2,
\end{equation}
where the discrete 2-norm of the 
numerical solution represents the energy.
The discrete version of
Parseval's relation relates the
2-norm of the numerical solution at a given time
to the 2-norm of its discrete Fourier coefficients by,
\begin{equation}
\|{\bf u}^{n}\|_2^2 =
\sum_{\mathbf{j}\in[N]^d }|u_\mathbf{j }^{n}|^2=\frac{1}{N^d}\sum_{\mathbf{m }\in[N]^d }|\widehat{u}_{\mathbf{m }}^{n}|^2=\frac{1}{N^d}\|\widehat{\bf u}^n\|_2^2.
\end{equation}
In the following two subsections
we consider options for the temporal discretization,
and
it show that each option yields a recurrence relation for the discrete
Fourier coefficients,
\begin{equation}
\widehat{u}_{\bf m}^{n+1}=\rho_{\mathbf m}\widehat{u}_{\bf m}^n,
\end{equation}
where $\rho_{\bf m}$ is the 
amplification factor of the
${\bf m}$-th discrete Fourier coefficient.
A sufficient condition for the
numerical method to satisfy
the discrete stability condition~\cref{eq:discrete_stability} is
\begin{equation} \label{eq:stab_cond}
|\rho_{\bf m}|\le1\iff |\widehat{u}_{\mathbf{m}}^{n+1}|\le |\widehat{u}_{\mathbf{
m}}^{n}| ~\text{for all wavenumbers}~\mathbf{m}\in[N]^d,
\end{equation}
which follows from application of Parseval's relation,
\begin{equation}
    \|{\bf u}^{n+1}\|_2^2
    =\frac{1}{N^d}\sum_{\mathbf{m }\in[N]^d }|\widehat{u}_{\mathbf{m }}^{n+1}|^2\le\frac{1}{N^d}\sum_{\mathbf{m }\in[N]^d }|\widehat{u}_{\mathbf{m }}^{n}|^2=\|{\bf u}^{n}\|_2^2.
\end{equation}

%%%%%%%%%%%%%%%%%%%

\subsection{Forward Euler scheme in time}

\begin{figure}[htb]
\centering
\includegraphics[width = .5\paperwidth]{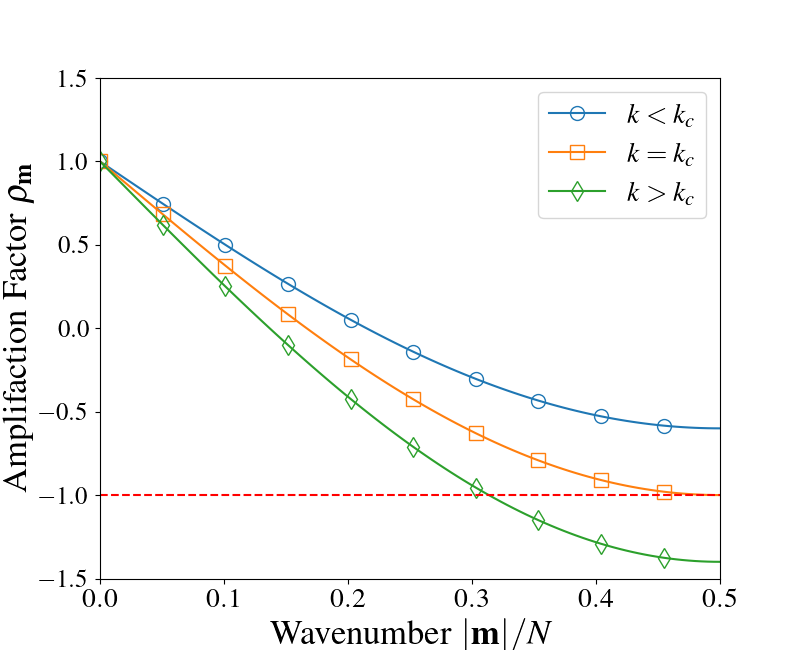}
\caption{
Forward Euler/central difference 
scheme~\cref{ds:HE_FD},
the amplification factor 
$\rho_{\bf m} = 1+kD\lambda_\mathbf{m}$ 
is plotted versus scaled wavenumber $|\mathbf{m}|/N$
for three time step values near the 
critical value,
$k < k_c, k = k_c, k > k_c$; 
for simplicity the spatial dimension is $d=1$.}
\label{fg:HE_amp}
\end{figure}

A forward difference in time yields the forward Euler/central difference scheme,
\begin{equation} \label{ds:HE_FD}
	\frac{u_\mathbf{j}^{n+1}-u_\mathbf{j}^n}{k}=D\nabla^2_hu_\mathbf{j}^n,
\end{equation}
where $k$ is the time step.
This can be expressed equivalently 
in terms of the discrete Fourier coefficients,
\begin{equation} \label{eq:HE_FD_FT}    \widehat{u}_{\mathbf{m}}^{n+1}=\widehat{u}_{\mathbf{m}}^{n}+k D\lambda_{\mathbf{m}}\widehat{u}_{\mathbf{m}}^{n} =
\rho_{\bf m}\widehat{u}_{\mathbf{m}}^{n},
\quad \rho_{\bf m} = 
(1+kD \lambda_{\mathbf{m}})~\text{for all}~\mathbf{m}\in[N]^d,
\end{equation} 
where $\rho_\mathbf{m}$ 
is the amplification factor of the scheme. 
The stability condition \cref{eq:stab_cond} then gives the following sufficient condition for stability of the scheme,
\begin{equation} \label{eq:stability}
    |1+k D\lambda_{\mathbf{m}}|\le1~ \text{for all}~\mathbf{m}\in[N]^d.
\end{equation}
Recalling the eigenvalue bound, $-4d/h^2\le\lambda_{\mathbf{m}}\le 0$,
in \cref{eq:apriori},
the stability condition~\eqref{eq:stability}
is guaranteed 
when the time step satisfies
\begin{equation}
|1+k D\lambda_{\mathbf{m}}| \le 1 
\Leftrightarrow
-1 \le 1+k D\lambda_{\mathbf{m}} \le 1 \Leftrightarrow
-2 \le k D\lambda_{\mathbf{m}} \le 0 \Leftrightarrow
k\le k_c,
\end{equation}
where $k_c=h^2/2dD$ is the critical time step.
Note that the lower bound in \cref{eq:apriori}
is achieved when $m_i=N/2, i=1,\ldots,d$
(assuming $N$ is even), 
so $k_c$ is a tight upper bound on stable time steps.
\Cref{fg:HE_amp} plots
the amplification factor of the scheme,
$\rho_{\bf m} = 1 + kD\lambda_\mathbf{m}$,
for three time step values,
$k < k_c, k = k_c, k > k_c$.
For simplicity the spatial dimension is $d=1$.
The results
indicate stability for $k\le k_c$
and
short wavelength instability for $k > k_c$.
For time step values $k > k_c$, 
the amplification factor lies outside
the stable range,
$\rho_{\bf m} < -1$,
for high wavenumber modes, 
indicating a short wavelength instability 
for the difference scheme in that case.
Next we consider a time discretization that does not
suffer from this instability.

%%%%%%%%%%%%%%%%%%%%%%%%%%%%%%%%%%%%%%%%%%%

\subsection{Crank-Nicolson scheme} 

\begin{figure}[htb]
\centering

\includegraphics[width = .5\paperwidth]{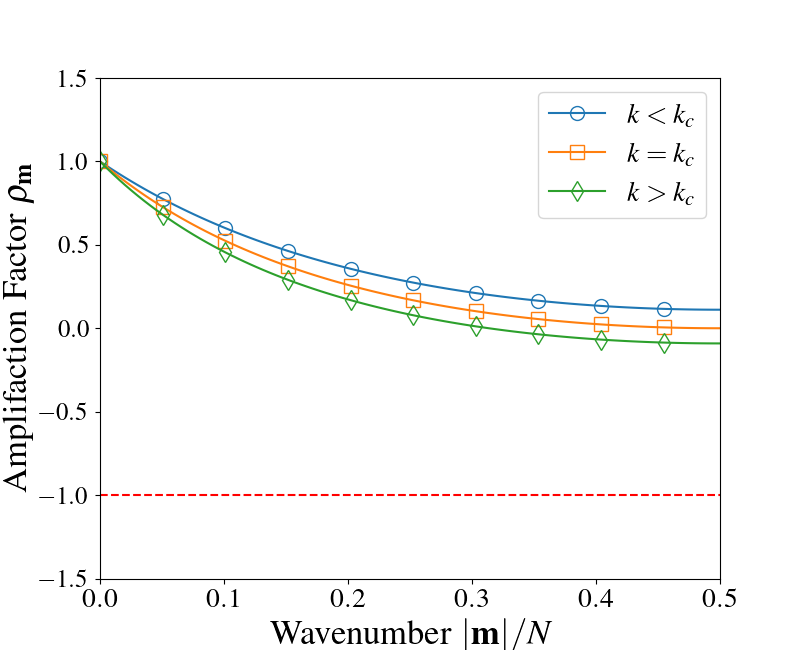}

\caption{
Crank-Nicolson scheme \cref{ds:HE_CN},
the amplification factor 
$\rho_{\bf m}
=(1+\frac{1}{2}kD\lambda_\mathbf{m})/(1-\frac{1}{2}kD\lambda_\mathbf{m})$ 
is plotted versus scaled wavenumber $|\mathbf{m}|/N$
is plotted for the same time steps $k$ as in \cref{fg:HE_amp}; for simplicity, the spatial dimension is $d=1$.
}
\label{fg:HE_amp_implicit}
\end{figure}

A trapezoidal difference in time yields the Crank-Nicolson scheme,
\begin{equation} \label{ds:HE_CN}
\frac{u_\mathbf{j}^{n+1}-u_\mathbf{j}^n}{k} =
\frac{1}{2}
D(\nabla^2_hu_\mathbf{j}^{n+1}+\nabla^2_hu_\mathbf{j}^n),
\end{equation}
which can be expressed equivalently in terms of the 
discrete Fourier coefficients of $u_{\bf j}^n$, 
\begin{equation}    \widehat{u}_{\mathbf{m}}^{n+1}=\widehat{u}_{\mathbf{m}}^{n}+ \frac{1}{2}kD(\lambda_{\mathbf{m}}\widehat{u}_{\mathbf{m}}^{n+1}+\lambda_{\mathbf{m}}\widehat{u}_{\mathbf{m}}^{n}),
\end{equation}
where the eigenvalues~$\lambda_{\mathbf{m}}$ of the 
discrete Lapacian were given in~ \cref{eq:eigenvalue}.
Then the amplification factor of the scheme is
obtained as follows,
\begin{equation} \label{ds:HE_CN_FT}
    \widehat{u}_{\mathbf{m}}^{n+1}= \frac{1+\frac{1}{2}kD\lambda_{\mathbf{m}}}{1- \frac{1}{2}kD\lambda_{\mathbf{m}}}\widehat{u}_{\mathbf{m}}^{n}=\rho_{\bf m}\widehat u_{\bf m}^n,
    \quad \rho_{\bf m}=\frac{1+\frac{1}{2}kD\lambda_{\mathbf{m}}}{1- \frac{1}{2}kD\lambda_{\mathbf{m}}}.
\end{equation}
Since $\lambda_{\mathbf{m}}\le0$, 
it follows that
$|1+\frac{1}{2}kD\lambda_{\mathbf{m}}| \le
|1-\frac{1}{2}kD\lambda_{\mathbf{m}}|$,
and then 
$|\rho_{\bf m}|\le1$ for all wavenumbers ${\bf m}$
and time steps $k > 0$,
so the scheme is unconditionally stable.
\Cref{fg:HE_amp_implicit} 
plots the amplification factor $\rho_{\bf m}$
of the 
Crank-Nicholson scheme for 
the same three time step values as in
\Cref{fg:HE_amp},
demonstrating the stability of
the scheme for these values.
Note also that
the Crank-Nicolson scheme \cref{ds:HE_CN} requires solving a linear system of $N^d$ equations 
at each time step.
However, even though \cref{ds:HE_CN} is an implicit equation for the numerical solution $u_{\bf j}^{n+1}$, 
\cref{ds:HE_CN_FT} 
gives an explicit expression for updating the discrete Fourier coefficients which scales like
$O(N^d)$,
and
transforming between the numerical solution
and 
its discrete Fourier coefficients
scales like $O(N^d\log N)$
using the Fast Fourier transform (FFT)~\cite{copetti_1990_kinetics}. 
Our solution of the linear system utilizes this FFT approach.

%%%%%%%%%%%%%%%%%%%%%%%%%%%%%%%%%%%

\subsection{Mass Conservation}

To check that both difference schemes satisfy 
the discrete analogue of mass conservation~\cref{eq:mass_con}, 
first note that 
the discrete Laplacian eigenvalue with wavenumber
$\mathbf{m = 0}$ is $\lambda_\mathbf{0}=0$.
Then the Fourier space equation~\cref{eq:HE_FD_FT} for the Forward Euler scheme becomes,
\begin{equation*}
\widehat{u}_\mathbf{0}^{n+1}=\widehat{u}_\mathbf{0}^n+kD\lambda_\mathbf{0}\widehat{u}_\mathbf{0}^n=\widehat{u}_\mathbf{0}^n,
\end{equation*}
so the average concentration is invariant in time,
\begin{equation} \label{eq:discrete_con}
	\frac{1}{N^d}\sum_{\mathbf{j}\in[N]^d}u_\mathbf{j}^{n+1}
    =\frac{1}{N^d}\widehat{u}_\mathbf{0}^{n+1}
    =\frac{1}{N^d}\widehat{u}_\mathbf{0}^{n}=\frac{1}{N^d}\sum_{\mathbf{j}\in[N]^d}u_\mathbf{j}^n.
\end{equation}
Similarly for the Crank-Nicolson scheme \cref{ds:HE_CN},
\begin{equation*}
	\widehat{u}_\mathbf{0}^{n+1}=\widehat{u}_\mathbf{0}^n+\frac{kD}{2}(\lambda_\mathbf{0}\widehat{u}_\mathbf{0}^{n+1}+\lambda_\mathbf{0}\widehat{u}_\mathbf{0}^n)=\widehat{u}_\mathbf{0}^n,
\end{equation*}
which also shows that the average concentration
is invariant in time.

%%%%%%%%%%%%%%%%%%%%%%%%%%%%%%%%%%%%%%

\section{Molecular Dynamics} \label{sec:MD}

The MD simulations presented below utilize
the Lennard-Jones potential to simulate the 
diffusion of argon in helium.
These are both noble gases, 
where argon has
atomic number $Z = 18$ and molar mass $m = 39.948$ g/mol,
while helium has 
atomic number $Z = 2$ and molar mass $m = 4.003$ g/mol.
This system represents one of the simplest real examples of chemical diffusion~\cite{Wasik1969-my}.
The short-range Lennard-Jones potential is
appropriate for modeling the diffusion of noble gases;
long-range electrostatic interactions
are more expensive to compute,
but they can be omitted 
in the case of small non-polar atoms such as these.
The experimentally determined value of the
diffusion coefficient $D_{\rm exp}$ of argon in helium
is available~\cite{Wasik1969-my}, which enables verification of the estimation algorithm developed in this work.
This section first discusses the equation governing
the system evolution,
then presents the Lennard-Jones potential,
and
finally gives the parameters used in the
MD simulations.

%%%%%%%%%%%%%%%%%%%%%%%%%

\subsection{Governing equation}

The argon and helium gases are represented by a set of
$N_p$ particles
$\{\mathbf{r}_i \in \mathbb{R}^2, i = 1,\ldots,N_p\}$
moving in two-dimensional space.
Newton's equations define the time evolution of the system,
\begin{equation} \label{eq:newton}
    m_i\frac{\ud^2\mathbf{r}_i}{\ud t}=\mathbf{F}_i,
    \quad i = 1, \ldots, N_p,
\end{equation}
where $m_i$ is the mass of particle $\mathbf{r}_i$ and $\mathbf{F}_i$ is the force on particle $\mathbf{r}_i$,
induced by all the other particles.
The total potential energy of the system
is the sum of  
pairwise particle interactions,
\begin{equation}
    V_\text{total}(\mathbf{r}_1, \ldots, \mathbf{r}_{N_p}) = 
    \sum_{i=1}^{N_p}\sum_{\substack{j = 1 \\ j\ne i}}^{N_p}V(r_{ij}), \quad
r_{ij}=|\mathbf{r}_i-\mathbf{r}_j|,
\end{equation}
where 
$V(r)$ is the
Lennard-Jones potential given below,
$r_{ij}$ is the 
interparticle distance between 
$\mathbf{r}_i$ and $\mathbf{r}_j$,
and the $j=i$ terms are excluded from the sum because particles do not interact with themselves.
The system is considered to be conservative
and
hence the force on each particle is the 
gradient of the total potential energy,
\begin{equation} 
\label{eq:MD_force}
\mathbf{F}_{i}=
\nabla_{\mathbf{r}_i}V_\text{total}(\mathbf{r}_1,\ldots,\mathbf{r}_{N_p})=
\sum_{i=1}^{N_p}\sum_{\substack{j = 1 \\ j\ne i}}^{N_p}
\nabla_{\mathbf{r}_{i}}V(r_{ij}),
\quad i = 1,\ldots,N_p,
\end{equation}
where $\nabla_{\mathbf{r}_{i}}$ is the gradient with respect to particle $\mathbf{r}_i$.

%%%%%%%%%%%%%%%%%%%%%%%

\subsection{Lennard-Jones potential}

\begin{figure}[h!]
\label{LJ_potential}
    \centering
    \subfigure[]{\includegraphics[width = .325\paperwidth]{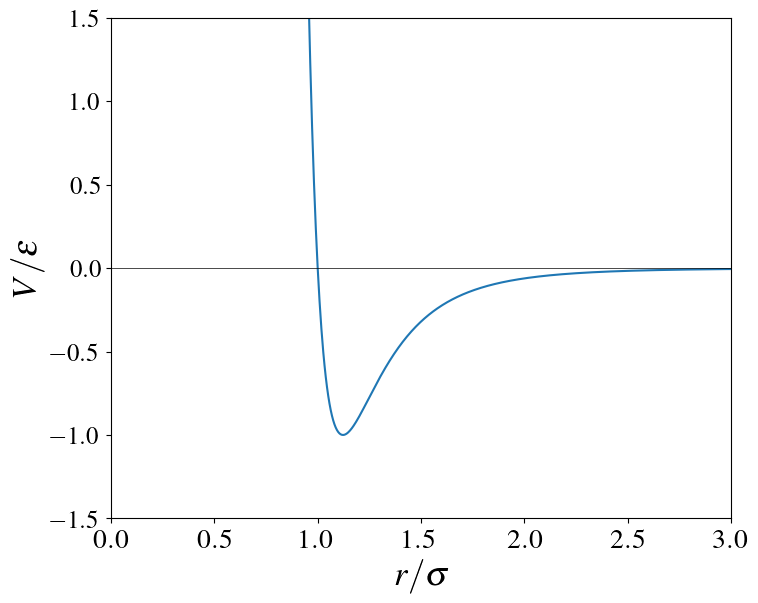}\label{sfg:lg_gen}} 
    \subfigure[]{\includegraphics[width = .325\paperwidth]{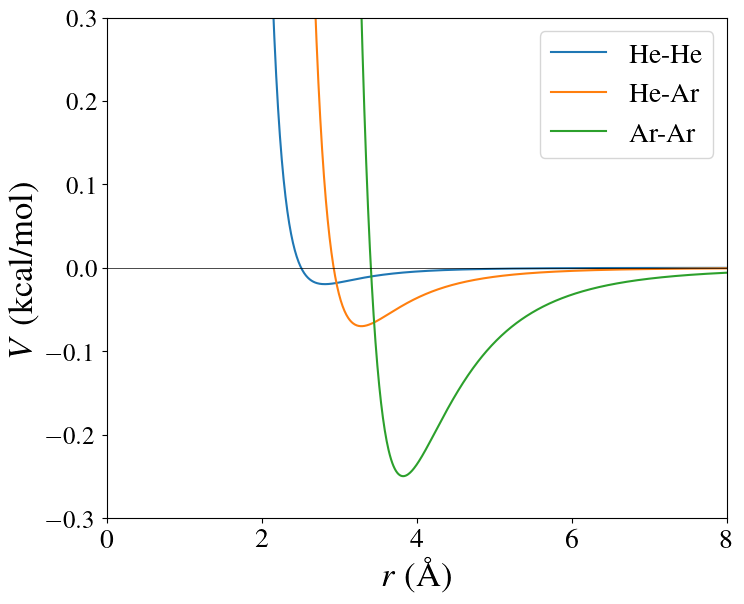}\label{sfg:lg_all}} 
\vskip -10pt
\caption{
Lennard-Jones potential $V$ given in \cref{eq:6-12_pot}, 
(a) scaled potential $V/\epsilon$
versus scaled distance $r/\sigma$,
(b) dimensional potentials $V$ (kcal/mol) for 
He-He, Ar-Ar, and mixed Ar-He interactions
versus interparticle distance $r$ (\r{A})
with 
parameters $\epsilon, \sigma$
given in \cref{tb:lj_diff}~\cite{helium,argon}.
}
\end{figure}

The Lennard-Jones potential is given by
the expression
\begin{equation} \label{eq:6-12_pot}
V(r)=
    \begin{cases}
        4\epsilon\left[ \left( \frac{\displaystyle\sigma}{\displaystyle r} \right)^{12}-\left( \frac{\displaystyle\sigma}{\displaystyle r} \right)^6 \right],&r<r_c,\\
        0,&r\ge r_c,
    \end{cases}
\end{equation}
where $\epsilon, \sigma$ are parameters, 
$r$ is the interparticle distance, 
and $r_c$ is a long-range cutoff. 
\Cref{sfg:lg_gen}
plots the 
scaled Lennard-Jones potential $V/\epsilon$
versus
the scaled interparticle distance
$r/\sigma$.
The potential characterizes 
London dispersion forces, where
the order 12~term
accounts for close-range repulsive forces
which dominate when $r<\sigma$,
and
the order 6~term accounts for
long-range attractive forces 
that operate for $r>\sigma$. 
The parameter $\epsilon$
sets the maximum depth of the potential well
which determines the 
interaction strength,
and
the parameter $\sigma$
sets the distance at which the 
potential $V(r)$ vanishes.

The Lennard-Jones parameters $\epsilon, \sigma$
depend on the physical and chemical properties 
of the particles in the system.
The parameters used in the present work
were taken from MD calculations
and
are displayed in \cref{tb:lj_diff}~\cite{helium, argon}.
In our simulations, 
each particle is designated as
either argon or helium,
and
the interaction potential $V(r)$ in \cref{eq:6-12_pot} between two particles is determined by their identities.
There are three types of pairwise interactions: helium/helium, argon/argon, 
and mixed helium/argon.
For each pair of particles,
the parameters are selected from~\cref{tb:lj_diff} and \cref{sfg:lg_all} plots the Lennard-Jones potentials for the three types of
pairwise interactions.
The relative distance scale of the helium/helium and argon/argon potentials reflects that argon has a larger atomic radius than helium.
Similarly, the argon/argon interaction is stronger and has a deeper well, because argon atoms are more polarizable, and exhibit greater induced dipoles.

\begin{table}[h!]
    \centering
      \begin{tabular}{|l|c|c|}
          \hline
          Interaction & $\epsilon$~(kcal/mol) & $\sigma$~(\r{A}) \\ \hline
          He-He       & 0.0196                     & 2.50                                   \\ \hline
          He-Ar       & 0.0700                     & 2.92                                   \\ \hline
          Ar-Ar       & 0.2498                     & 3.40                                   \\ \hline
      \end{tabular}  
      \vskip 2.5pt
\caption{Lennard-Jones parameters for argon and helium as determined in 
MD calculations~\cite{helium, argon},
mixed parameter values use
the geometric mean of the
individual particle parameters.
The cutoff distance is $r_c =20$\r{A}
for all pairwise interactions.}
\label{tb:lj_diff}
\end{table}

%%%%%%%%%%%%%%%%%%%%%

\subsection{LAMMPS} 

This work utilizes the 
LAMMPS software package
to implement the MD simulations
for argon diffusion in helium
with the
Lennard-Jones potentials described above~\cite{LAMMPS}.
The real unit system in LAMMPS
is chosen for these simulations,
where
distance is in Angstroms (\r{A}), 
energy in kilocalories per mole (kcal/mol), 
time in femtoseconds (fs = 1e-15 second), 
temperature in Kelvin (K), 
and
mass in grams per mole (g/mol).
The simulations were done in
two space dimensions,
where the MD simulation box 
has dimensions 
5e4 \r{A} $\times$ 5e4 \r{A},
but for simplicity 
the results are presented using
a dimensionless length scale where 
the unit square $[0,1]^2$
represents this domain.
Periodic boundary conditions were imposed
using the minimum image convention
in which the particles are replicated across the
sides of the simulation box for the purpose of
computing the forces, 
and upon exiting a side, 
they re-enter the opposite side.
Three MD simulations were performed
in which the initial condition in each 
is a unique pseudo-random seed 
of 3e4 helium atoms uniformly
distributed in the domain $[0,1]^2$
and a patch of 3e4 argon atoms uniformly
distributed in the square subdomain $[1/4,3/4]^2$.
The initial particle velocities were chosen from the Maxwell-Boltzmann distribution,
\begin{equation}
f(\mathbf{v}) \propto
\e^{\displaystyle -\frac{m|\mathbf{v}|^2}{2k_bT}},
\end{equation}
where $T=300$ K, $m$ is the mass of a particular atom, and $k_b=0.001987$ kcal/(mol K) is Boltzmann's constant. 
The MD timestep was $k_\text{MD}=5$ fs. 
The order 12 repulsive term of the Lennard-Jones potential \cref{eq:6-12_pot} necessitates a small timestep 
to ensure sufficient accuracy in the 
LAMMPS integrator. 
An NVE integrator was used for
the MD simulations, which samples the microcanonical ensemble, where the number of particles, volume, and total energy are fixed.
The LAMMPS NVE integrator uses the
St\"ormer-Verlet time stepping scheme.

\section{Comparison of MD and FD results} \label{sec:res}

This section compares molecular dynamics
and 
finite-difference
simulations of argon diffusion in helium
in two space dimensions.
Recalling that the MD simulations utilized
random initial conditions representing a
square patch of argon atoms,
the FD calculations were initialized
with the following continuum idealization,
\begin{equation}
\label{eq:FD_initial}
u_0(\mathbf{x})=\begin{cases} 1, & \mathbf{x}\in[1/4,3/4]^2, \\
0, & \mathbf{x}\in[0,1]^2\setminus[1/4,3/4]^2.
\end{cases}
\end{equation}
The total simulation time was $5~{\rm ns}$,
where
the FD time step was 
$k_{\rm FD} = 5~{\rm ps}$ (1e3 time steps)
and
the MD time step was
$k_{\rm MD} = 5~{\rm fs}$ (1e6 time steps).
The MD time step $k_{\rm MD}$ 
was determined empirically to achieve stable results
and is 1000 times smaller
than the FD time step $k_{\rm FD}$;
this is because 
the MD simulations with the explicit
St\"ormer-Verlet time stepping scheme
are constrained
by the stiffness of the 
Lennard-Jones potential,
whereas the FD calculations
use the implicit Crank-Nicolson time stepping scheme
which is unconditionally stable.
To accommodate the difference between the MD and FD time steps,
the fitting procedure discussed in the next section
samples the MD trajectory every 1000 time steps.

The simulations were carried out on 
a single CPU core of
an Apple M1 2020 Macbook Pro with 16 GB of memory running the Sonoma 14 operating system.
The run time for each FD calculation was less than a second,
and
the run time for each MD simulation 
was about 80 minutes.
A parallel implementation
would reduce the MD run time significantly,
and
LAMMPS does support parallel computing
on multicore systems,
but this is beyond the scope of the present work.

\cref{fg:diff_results}
compares the atomistic MD results
and
the continuum FD results computed with grid spacing $h = 1/50$.
The diffusion coefficient
in the FD simulations ($D = 0.7948$~cm$^2$/s)
is the optimal value
determined by the estimation algorithm
described in the next section.
The time is given at the top of each column.
Row 1 plots the MD results,
where the argon atoms are red dots
and 
the helium atoms are a uniform blue background.
As time proceeds,
the initial square distribution of argon atoms
evolves into an expanding radially symmetric 
distribution.

\def\heightpad{31pt}
\begin{figure}[htbp]

    \centering
    \begingroup
    \setlength{\tabcolsep}{1pt}
    \def\subheight{.117\paperwidth}
    \begin{tabular}{rccccc}

        ~ & $t=0$ ns & $t=0.25$ ns & $t=1$ ns & $t=2.5$ ns & $t=5$ ns \\
        
        1 MD &
        \includegraphics[align = c, width=\subheight]{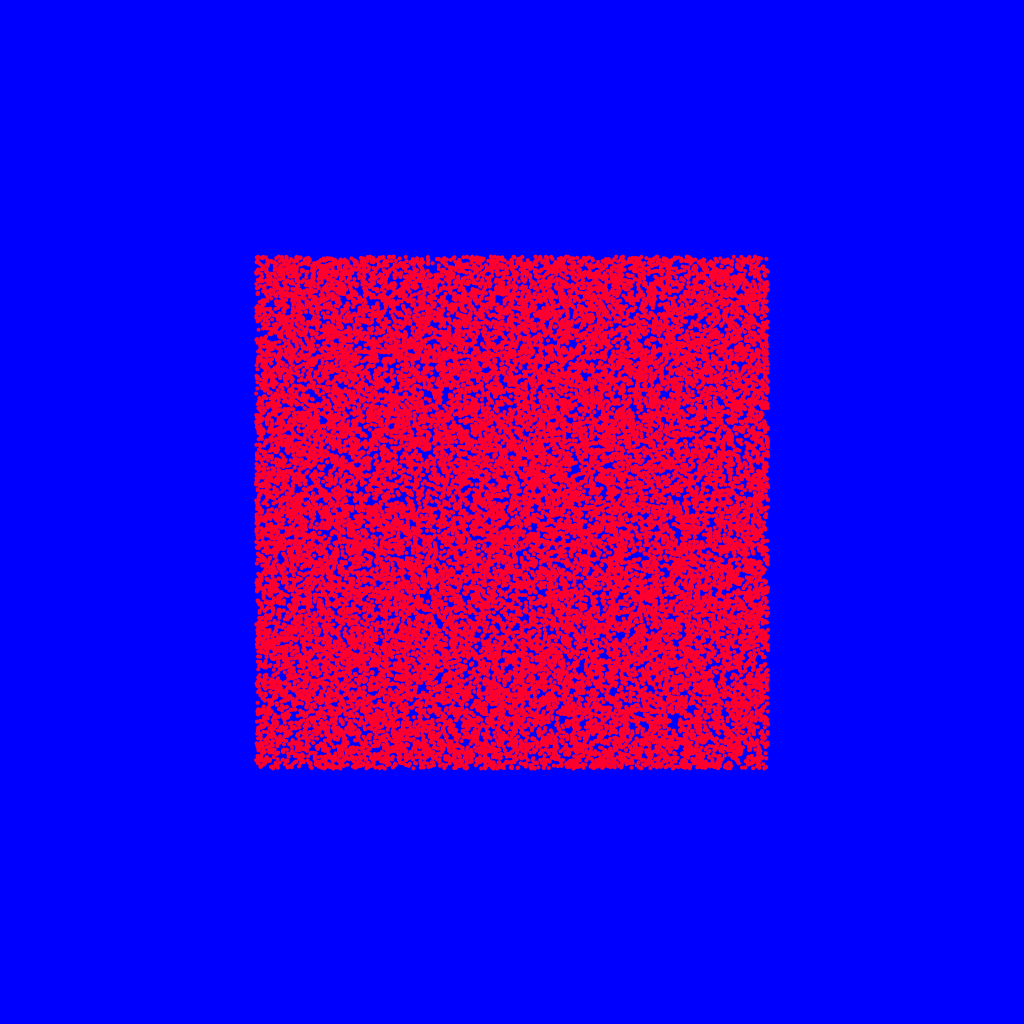} & 
        \includegraphics[align = c, width=\subheight]{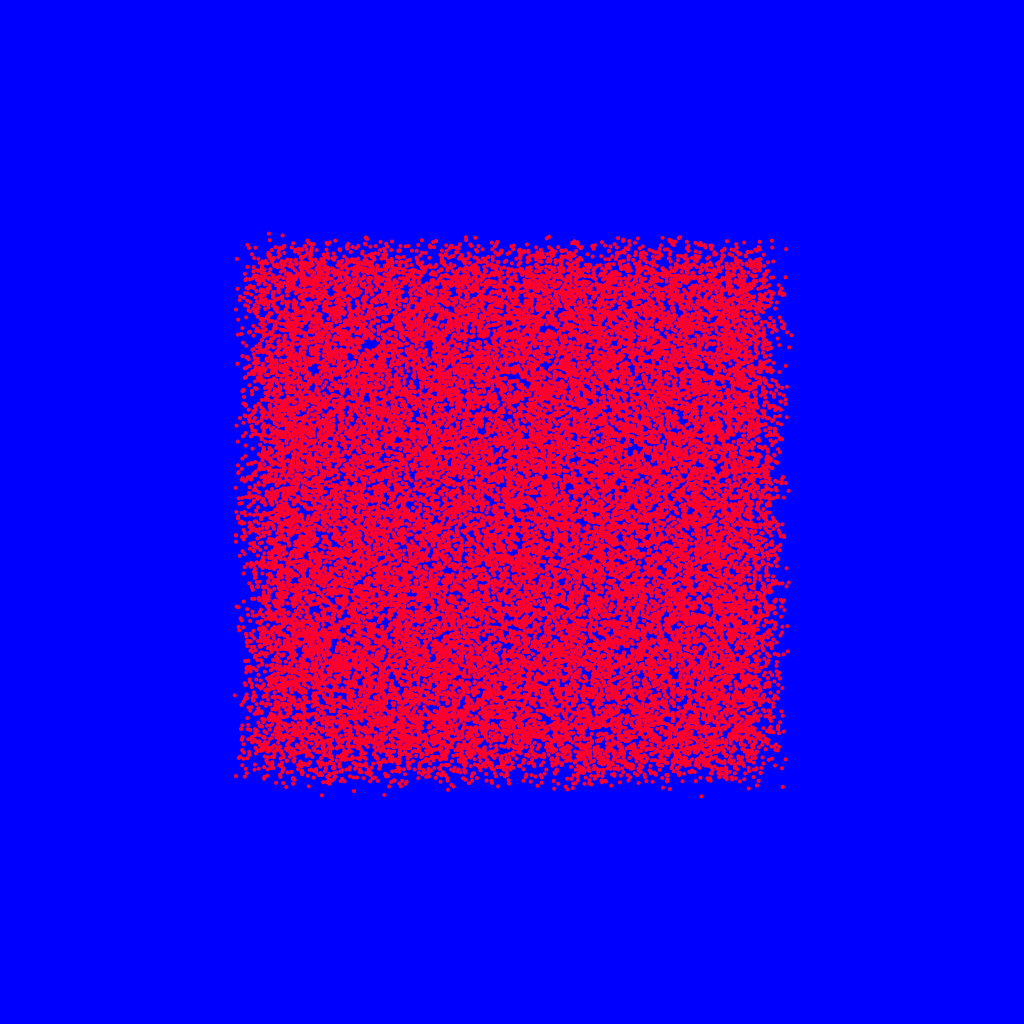} & 
        \includegraphics[align = c, width=\subheight]{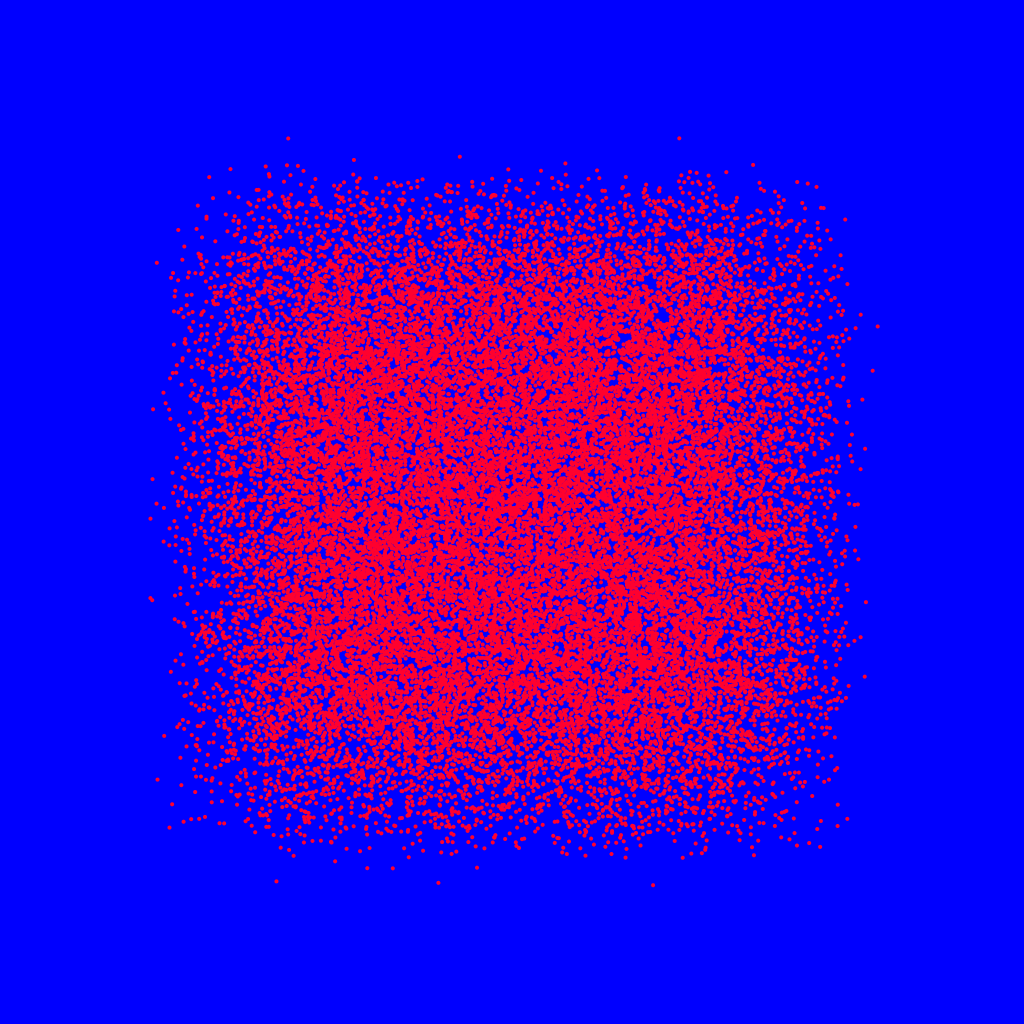} & 
        \includegraphics[align = c, width=\subheight]{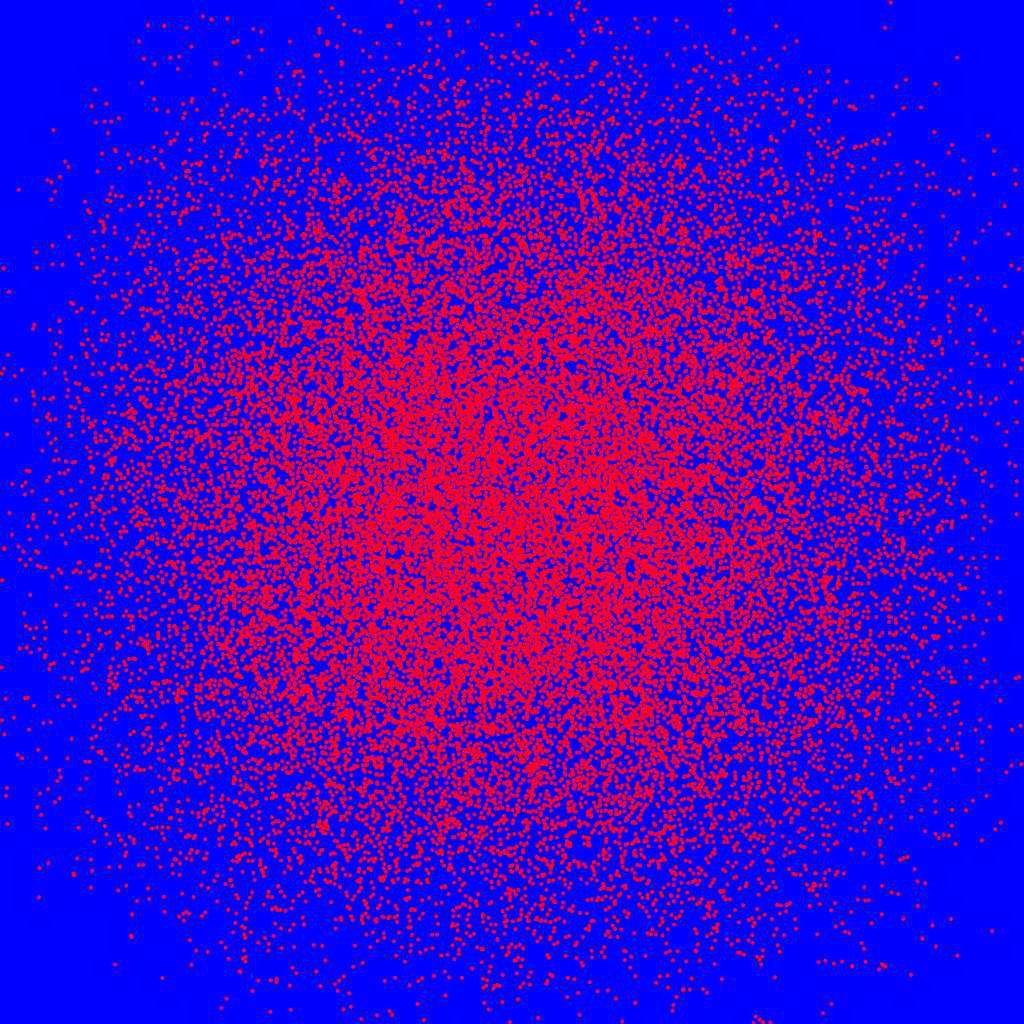} & 
        \includegraphics[align = c, width=\subheight]{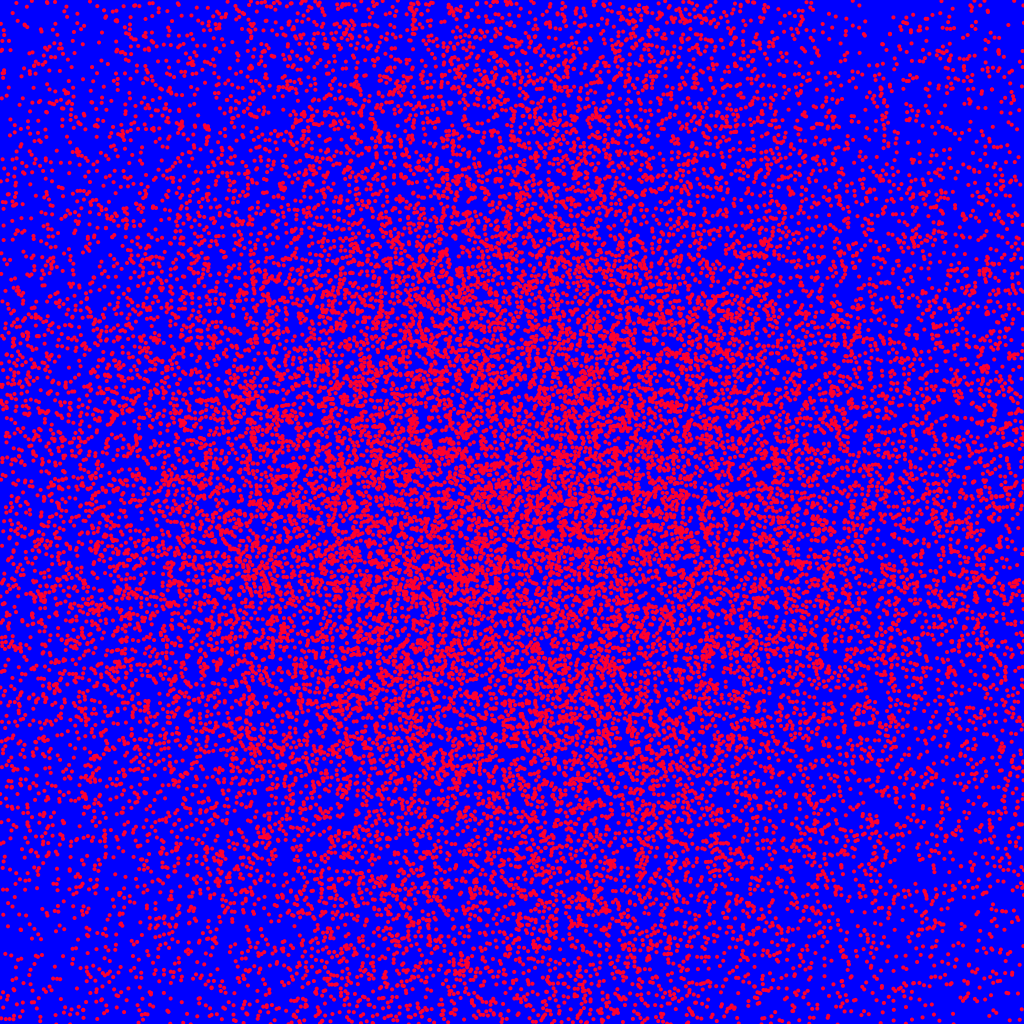} \\[\heightpad]

        2 MD & 
        \includegraphics[align = c, height=\subheight]{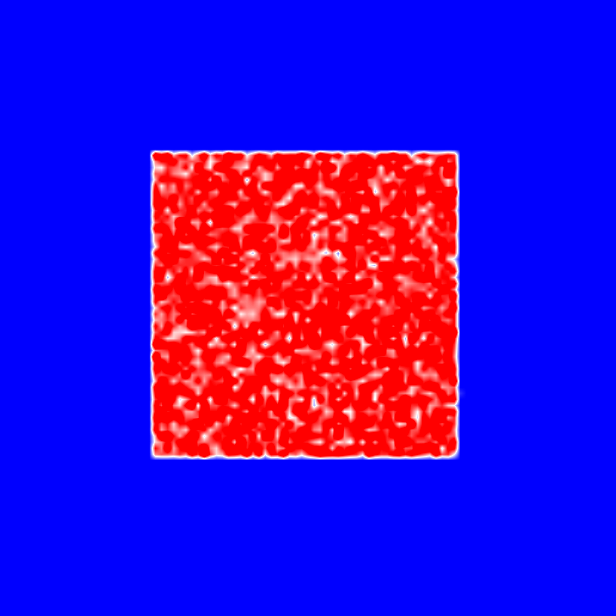} & 
        \includegraphics[align = c, height=\subheight]{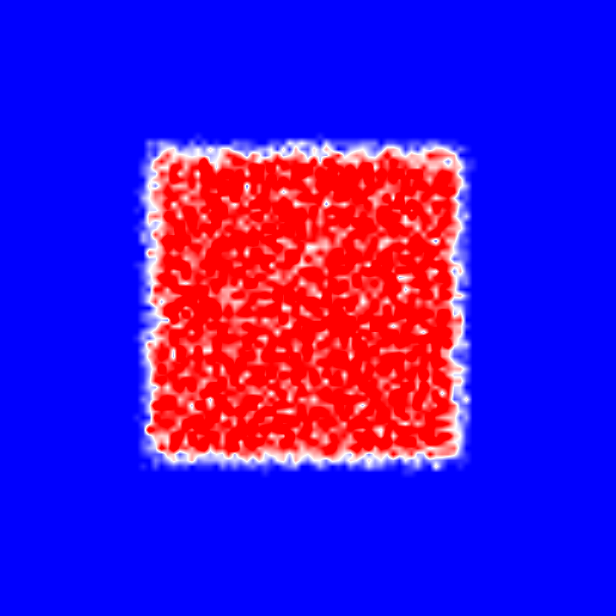} & 
        \includegraphics[align = c, height=\subheight]{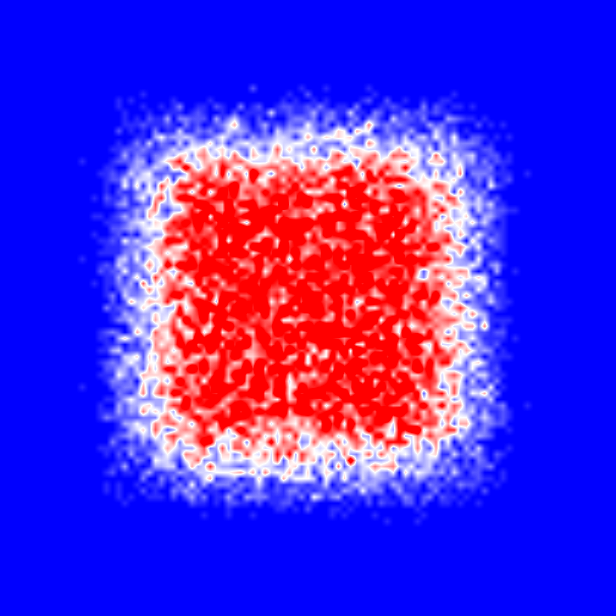} & 
        \includegraphics[align = c, height=\subheight]{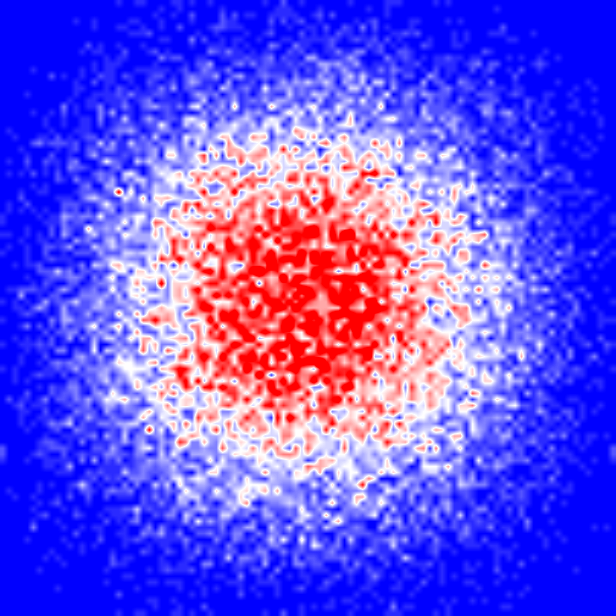} & 
        \includegraphics[align = c, height=\subheight]{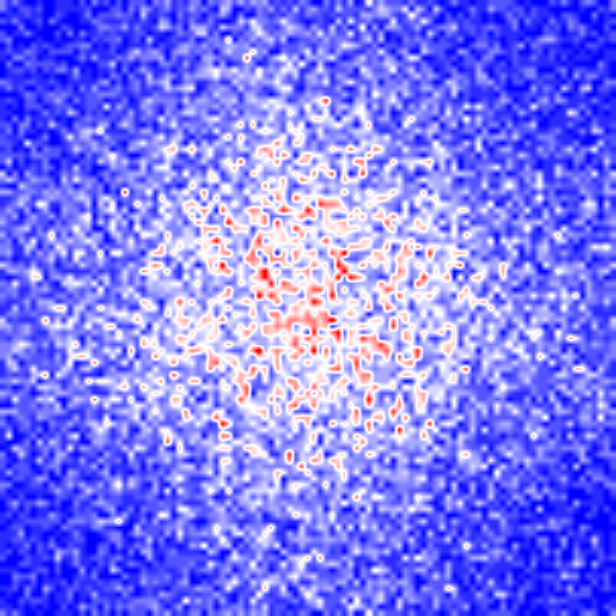} \\[\heightpad]

        3 FD &
        \includegraphics[align = c, height=\subheight]{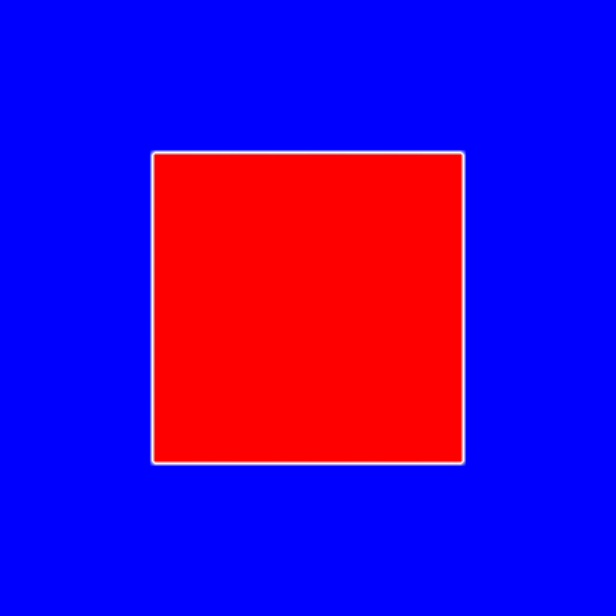} & 
        \includegraphics[align = c, height=\subheight]{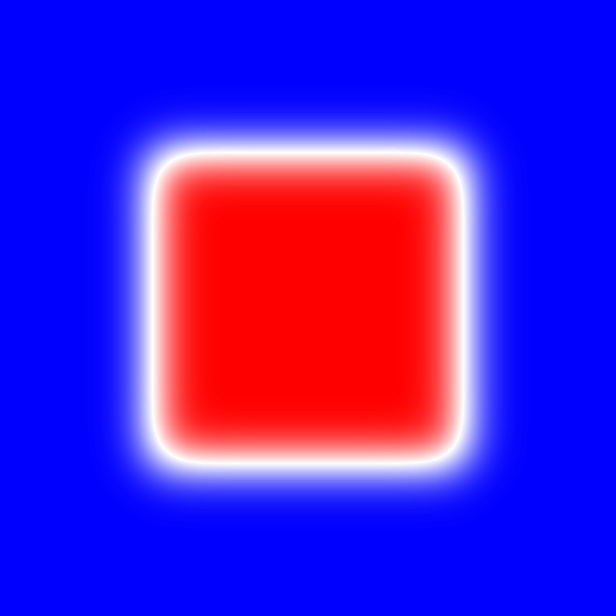} & 
        \includegraphics[align = c, height=\subheight]{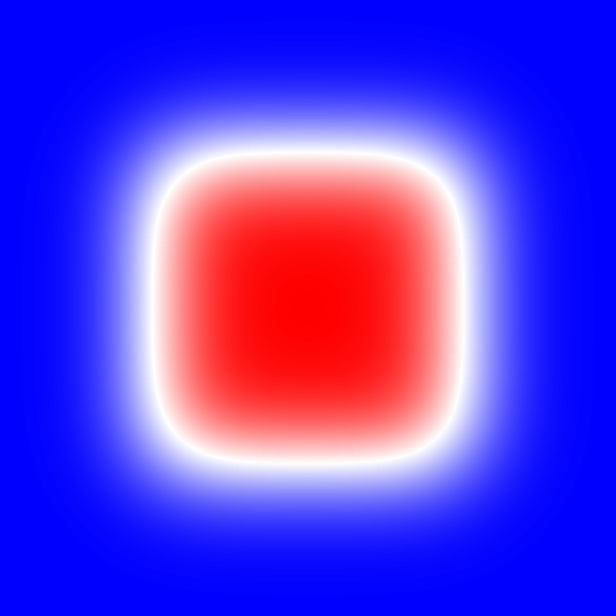} & 
        \includegraphics[align = c, height=\subheight]{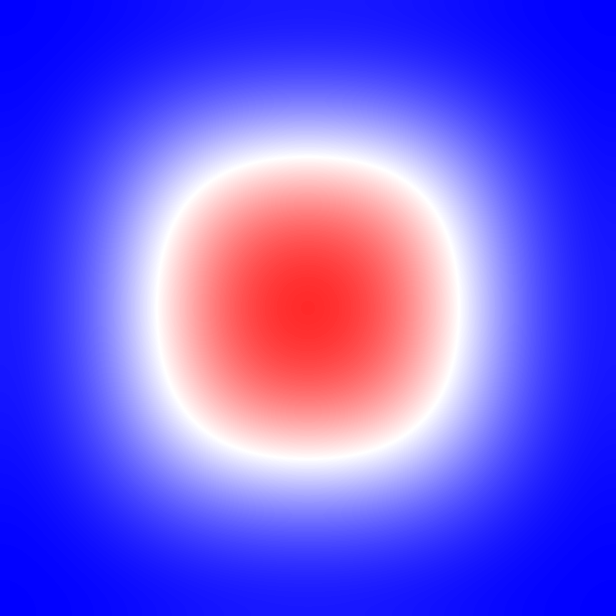} & 
        \includegraphics[align = c, height=\subheight]{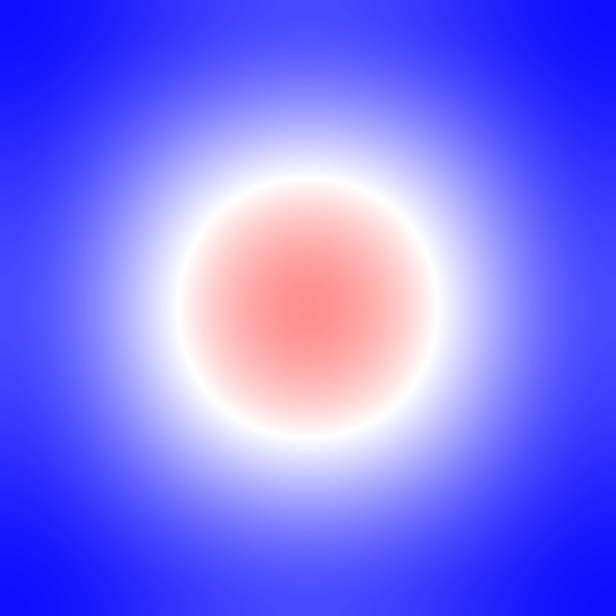} \\[\heightpad]

        4 MD &
        \includegraphics[align = c, height=\subheight]{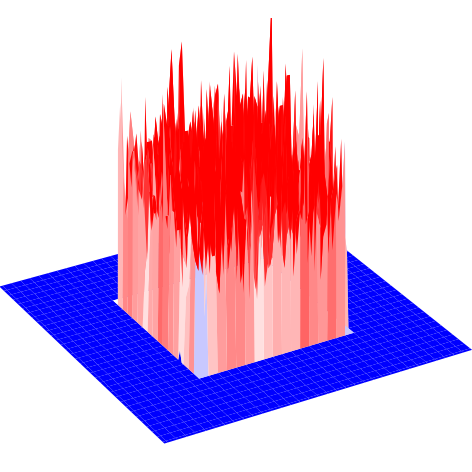} &
        \includegraphics[align = c, height=\subheight]{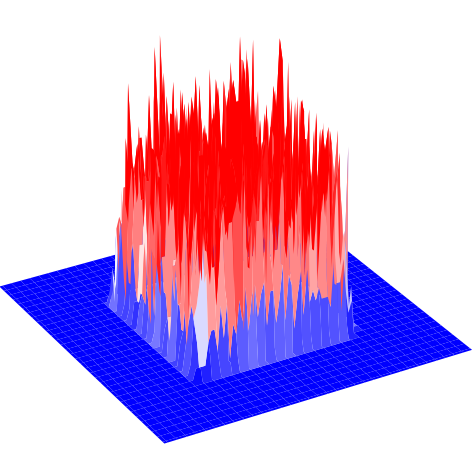} &
        \includegraphics[align = c, height=\subheight]{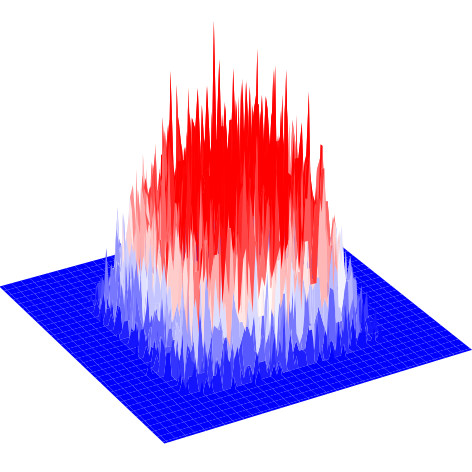} &
        \includegraphics[align = c, height=\subheight]{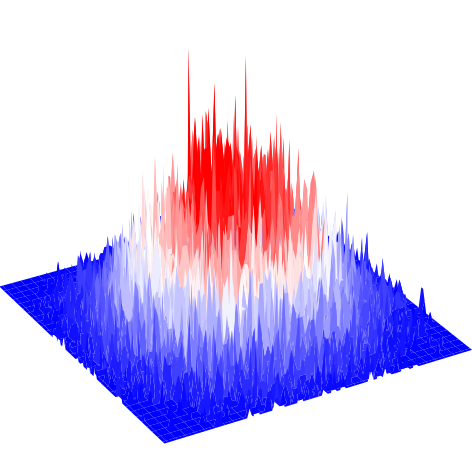} &
        \includegraphics[align = c, height=\subheight]{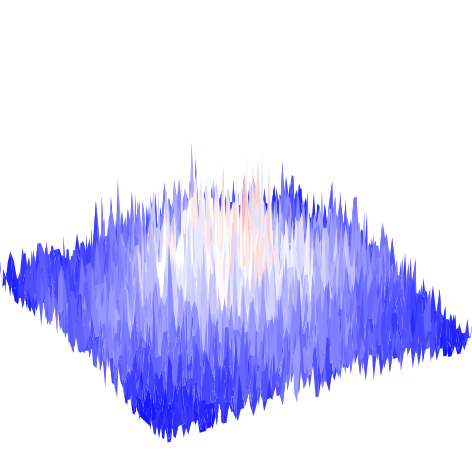} \\[\heightpad]

        5 FD &
        \includegraphics[align = c, height=\subheight]{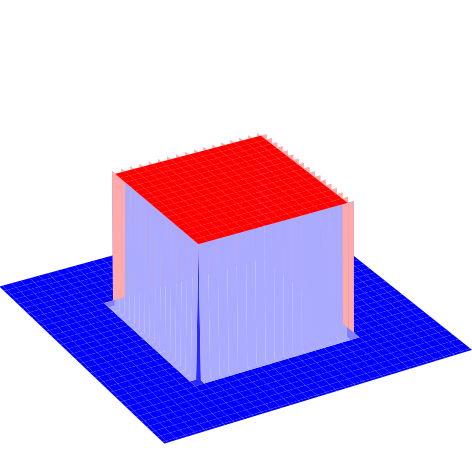} &
        \includegraphics[align = c, height=\subheight]{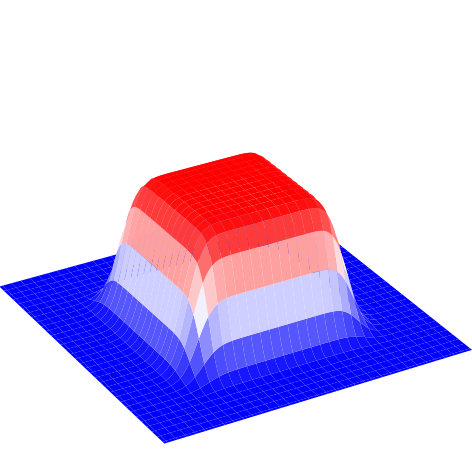} &
        \includegraphics[align = c, height=\subheight]{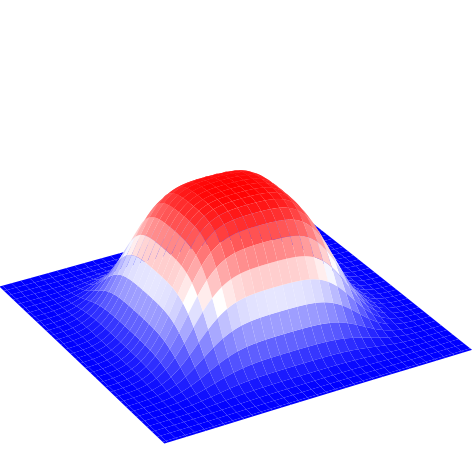} &
        \includegraphics[align = c, height=\subheight]{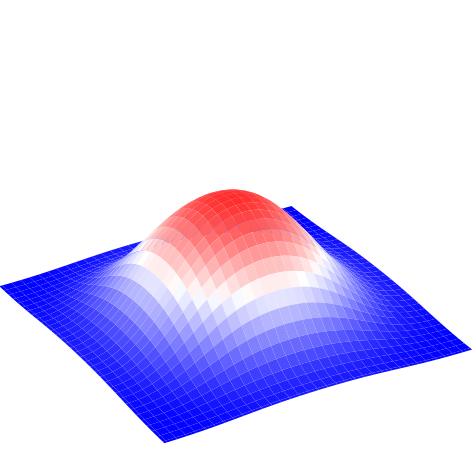} &
        \includegraphics[align = c, height=\subheight]{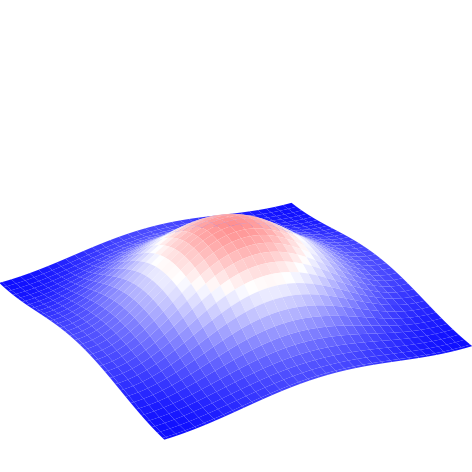} \\[\heightpad]

        6 MD &
        \includegraphics[align = c, width=\subheight]{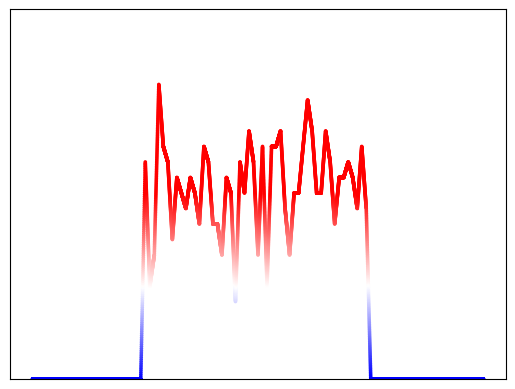} &
        \includegraphics[align = c, width=\subheight]{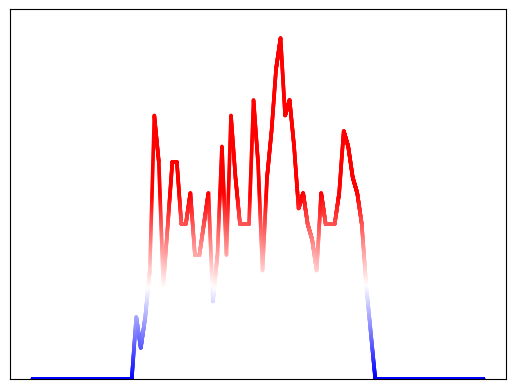} &
        \includegraphics[align = c, width=\subheight]{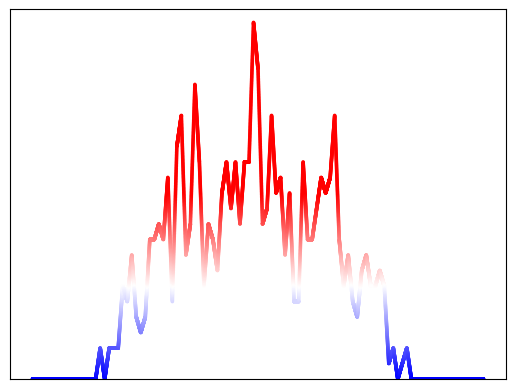} &
        \includegraphics[align = c, width=\subheight]{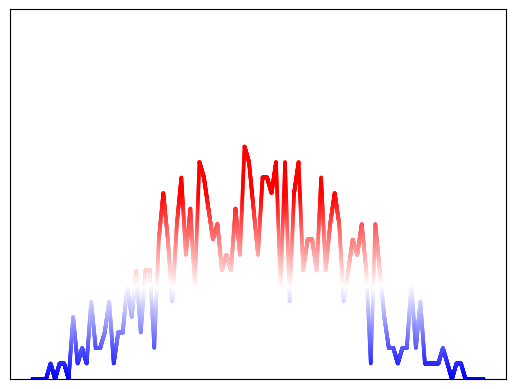} &
        \includegraphics[align = c, width=\subheight]{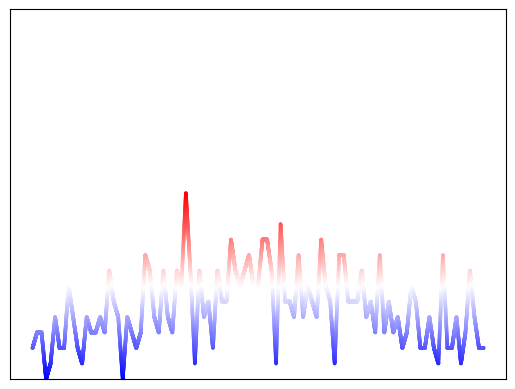} \\[\heightpad]

        7 FD &
        \includegraphics[align = c, width=\subheight]{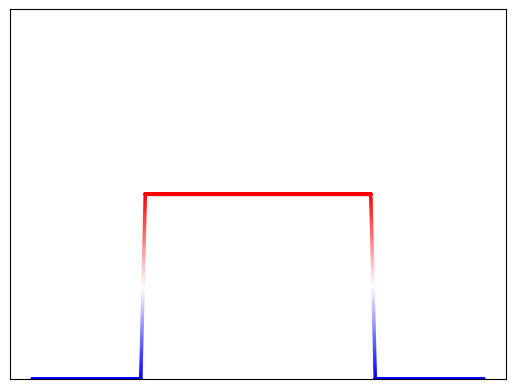} &
        \includegraphics[align = c, width=\subheight]{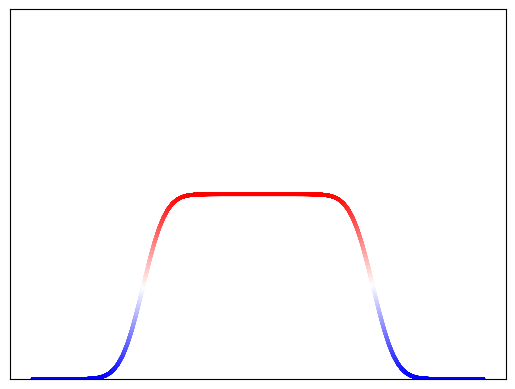} &
        \includegraphics[align = c, width=\subheight]{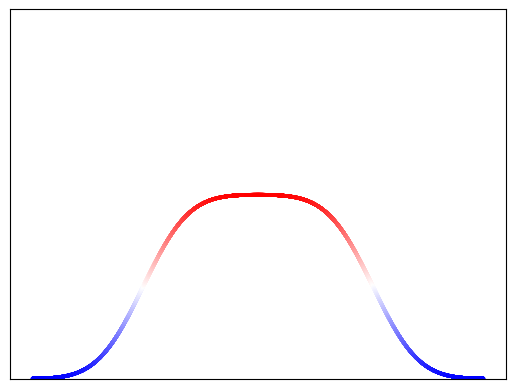} &
        \includegraphics[align = c, width=\subheight]{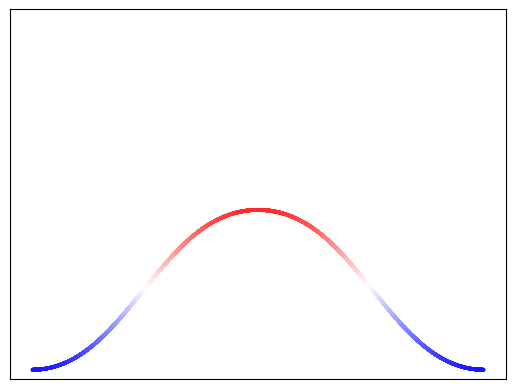} &
        \includegraphics[align = c, width=\subheight]{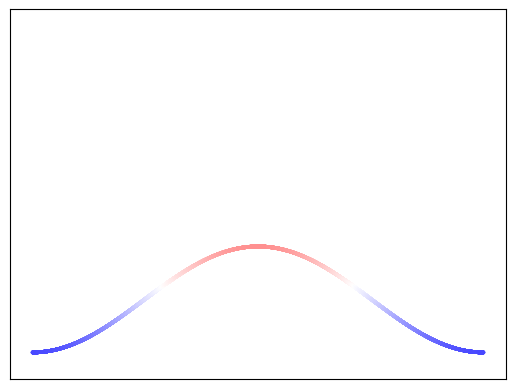}
    \end{tabular}
    \endgroup

    \hspace{25pt}
    \includegraphics[width = .53\paperwidth]{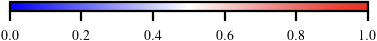}

    \caption{
Numerical results for diffusion of
argon in helium,
color bar indicates normalized argon concentration,
MD binning parameter $N=50$,
FD solution of diffusion equation~\cref{eq:HE}
computed on a $50 \times 50$ grid,
Columns 1-5 display times $t=0,0.25,1,2.5,5$ ns,
Row 1: MD snapshots of argon atoms,
Rows 2, 4, 6: binned MD concentration,
Rows 3, 5, 7: computed FD concentration.
}
\label{fg:diff_results}
\end{figure}

To facilitate comparison of the MD and FD results,
the MD argon trajectories $\{{\bf r}_i^n\}$
were converted into a concentration $U_{\bf j}^n$ 
defined on the FD grid
by
a binning process~\cite{larson_1997_hydrodynamics}
that counts the 
number of argon atoms in each FD grid cell,
$N_{\bf j}^n$,
and
scales the count to lie between 0 and 1,
\begin{equation} \label{eq:bin}
N_\mathbf{j}^n = 
\#\{\mathbf{r}_i^n\in \text{grid cell }\mathbf{j}\},
\quad
U_{\mathbf{j}}^n=\frac{N_\mathbf{j}^n}{\max N_\mathbf{j}^n}.
\end{equation} 
Row 2 in~\cref{fg:diff_results} plots the 
resulting binned MD argon concentration $U_{\bf j}^n$
and
Row 3 plots the 
corresponding FD concentration $u_{\bf j}^n$
as heat maps on the FD grid.
Although the MD results are noisy
and
the FD results are smooth,
there is a clear similarity between the two
concentrations.
The similarity is also seen in
Rows 4 and 5 which present surface plots of the MD and FD concentrations.
Rows 6 and 7 show the $x = 1/2$ cross-section of the MD and FD concentrations.
In the surface plots and cross-sections, 
the MD and FD maximum argon concentrations 
decrease at similar rates.
\cref{fg:diff_results}
shows that using an optimal value of the 
diffusion coefficient $D$ in the FD calculation
ensures a close connection between the FD and MD results.

One may consider how the match between
the MD and FD results depends
on the binning parameter $N$,
which also determines the FD grid size,
$h = 1/N$.
The accuracy of the FD solution $u_{\rm j}^n$
relative to the true solution 
of the diffusion equation $u({\bf x},t)$
increases with $N$.
Additionally, 
the binned MD concentration $U_{\bf j}^n$
will better resolve the computed MD trajectories ${\bf r}_i^n$ for large $N$.
However as the binning parameter $N$ increases, 
the binned MD concentration will become
increasingly noisy
due to the use of a finite number particles in
the MD simulation.
Based on these competing effects, we expect there to be an optimal level of binning resolution to achieve the best match between the
MD and FD results.
The balance between increasing resolution and noise also depends on several other factors
including the MD domain size, 
number of MD particles, 
and initial MD distribution.

%%%%%%%%%%%%%%%%%%%%%%%%%%

\section{Diffusion coefficient 
estimation algorithm} \label{sec:estim}

This section describes the algorithm for
estimating the optimal diffusion coefficient
$D_{\rm opt}$
for the continuum diffusion equation~\cref{eq:HE}.
This is essentially a 
least-squares fitting procedure that
extracts $D_{\rm opt}$ from the
MD trajectories.
Recall that
the MD trajectories~${\bf r}_i^n$
give the position of each 
argon and helium atom at a 
sequence of discrete time steps,
and
these trajectories are represented on the FD grid
by the binned MD argon concentration $U_\mathbf{j}^n$
obtained using the previously described binning process.
Let $u_\mathbf{j}^n(D)$
be the numerical solution of the
continuum diffusion equation~\cref{eq:HE}
with diffusion coefficient $D$
given by the Crank-Nicolson scheme~\cref{ds:HE_CN}.
The vector of residuals is given by
\begin{equation} \label{eq:residuals}
\{ U_\mathbf{j}^n - u_\mathbf{j}^n(D) \}, \quad
n=0,\ldots,n_\text{max}, \quad
\mathbf{j} \in [N]^2,
\end{equation}
and the cost function is the
mean square deviation between
the MD and FD results,
\begin{equation} \label{eq:cost}
\texttt{cost}(D) =
\frac{1}{N^2}\sum_{n=0}^{n_\text{max}}\sum_{\mathbf{j}\in [N]^2}(U_{\mathbf{j}}^n - u_{\mathbf{j}}^n(D))^2.
\end{equation}
Minimizing the cost function
is a nonlinear optimization problem 
that yields the optimal diffusion coefficient,
$D_\text{opt}(N)$.
\Cref{fg:cost} plots the cost function
for the argon-helium system with 
MD binning parameter $N=20, 50, 100$,
showing that $\texttt{cost}(D)$ is a smooth convex function with a 
unique global minimum $D_\text{opt}(N)$
in each case.

\newcommand{\costwidth}{.3\textwidth}
\begin{figure}[h]
\centering
\begin{tabular}{ccc}
\includegraphics[width = \costwidth]{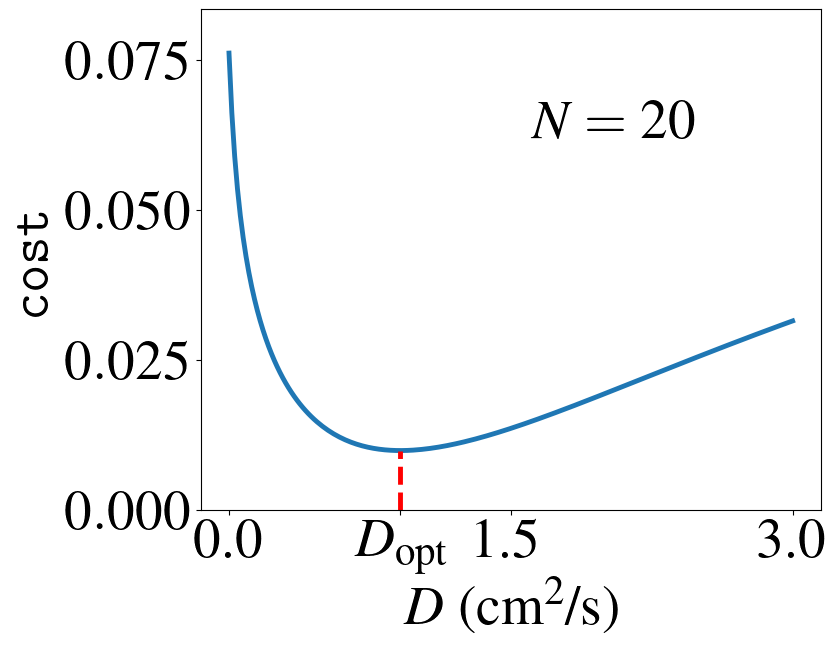} & 
\includegraphics[width = \costwidth]{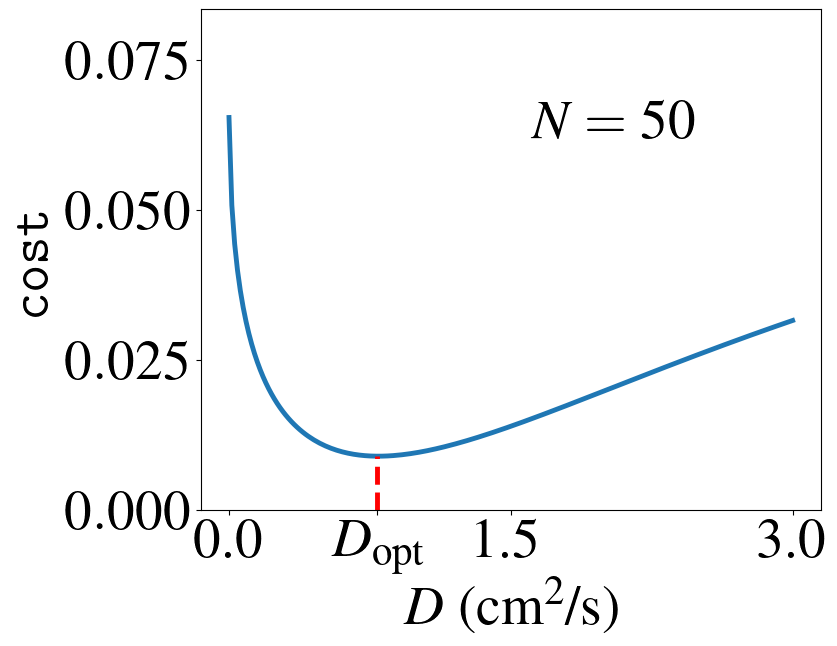} & 
\includegraphics[width = \costwidth]{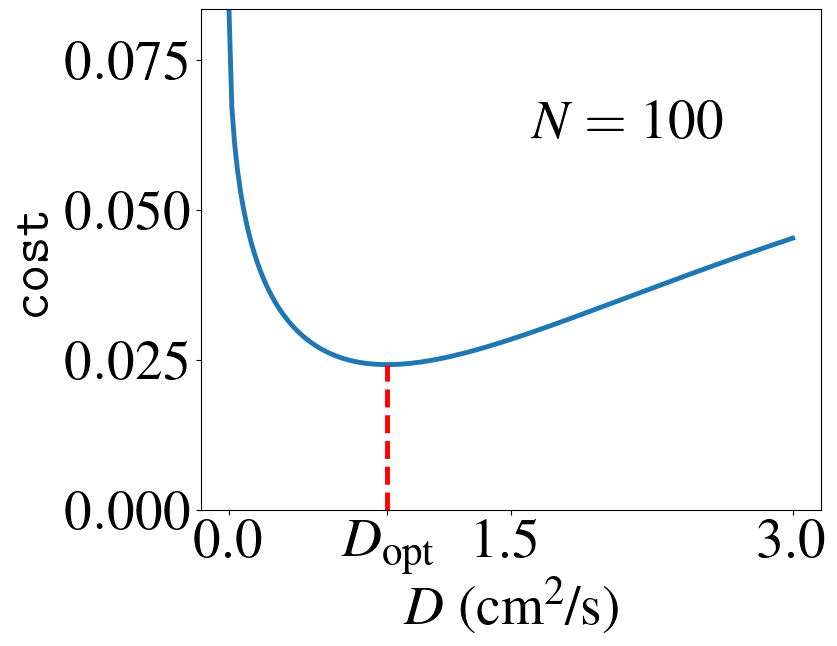} 
\end{tabular}
\caption{The cost function 
${\tt cost}(D)$~\cref{eq:cost}
of the diffusion coefficient 
estimation algorithm
is plotted versus the diffusion coefficient $D$
near the global minimum $D_\text{opt}(N)$ 
with MD binning parameters $N=20, 50, 100$.}
\label{fg:cost}
\end{figure}

\Cref{alg:estimation} shows
the steps for computing the
optimal diffusion coefficient,
\begin{equation} \label{eq:lstq}
D_{\rm opt}(N) = \arg\min_D \texttt{cost}(D) =
\arg\min_D\frac{1}{N^2}\sum_{n=0}^{n_\text{max}}\sum_{\mathbf{j}\in [N]^d}(U_{\mathbf{j}}^n - u_{\mathbf{j}}^n(D))^2.
\end{equation}
Lines 1-3 input the 
MD argon trajectories $\mathbf{r}_i^n$,
MD binning parameter $N$,
and 
an initial guess for the
diffusion coefficient $D_0$.
Line 4 computes the binned MD concentration,
$U_\mathbf{j}^n$,
using the binning procedure described in \cref{eq:bin}.
Line 5 computes the 
optimal diffusion coefficient $D_{\rm opt}(N)$
using the \texttt{curve\_fit} routine
from the Python package \texttt{scipy}~\cite{scipy}
to solve the 
least-squares minimization problem,
Line 6 returns $D_{\rm opt}(N)$.
Next we explain how the
\texttt{curve\_fit} routine works.
 
\begin{algorithm}[htb] 
\caption{Diffusion Coefficient Estimation Algorithm}
\label{alg:estimation}
\begin{algorithmic}[1]
\State \textbf{input} MD trajectories, 
$\mathbf{r}_i^n, i=1,\ldots,N_p$,
time steps $n=0,\ldots,n_\text{max}$
\State \textbf{input} MD binning parameter, $N$
\State \textbf{input} initial diffusion coefficient guess, $D_0$
\State compute binned MD concentration,
$U_\mathbf{j}^n$, from \cref{eq:bin}
\State compute $D_{\rm opt}(N)$ from \cref{eq:lstq}
using the \texttt{curve\_fit}
routine from \texttt{scipy}
\State \Return $D_{\rm opt}(N)$
\end{algorithmic}
\end{algorithm}

The $\texttt{curve\_fit}$ routine uses the Levenberg-Marquardt (LM) algorithm~\cite{LM} to compute the optimal diffusion 
coefficient~$D_{\rm opt}(N)$. 
The input for $\texttt{curve\_fit}$ 
is an initial guess $D_0$
for the diffusion coefficient.
The LM algorithm is an iterative scheme
that combines the stability of
gradient descent
for minimizing the 
cost function $\texttt{cost}(D)$
with the accuracy of Newton's method
for finding a
critical point of $\texttt{cost}(D)$.
At each step,
the diffusion coefficient is updated,
$D \to D + \delta D$,
where the increment $\delta D$
is the solution of the normal equations,
\begin{equation}
\label{eq:LM}
(\mathbf{J}^\text{T}\mathbf{J} + \lambda)
\delta D = \mathbf{J}^{\rm T}(U_\mathbf{j}^n-u_\mathbf{j}^n(D)),
\end{equation}
that arise from linearizing
the least-squares minimization problem~\cref{eq:lstq}.
In these equations,
$\mathbf{J} = \partial u_\mathbf{j}^n/\partial D$ 
is the Jacobian of the residual vector,
which is computed numerically,
and
$\lambda$ is an adjustable
regularization parameter. 
When the current step causes a 
large reduction in ${\tt cost}(D)$, 
$\lambda$ is set to a small value
and 
the algorithm is close to Newton's method,
but if the reduction in ${\tt cost}(D)$ is small,
then $\lambda$ is increased to bring the 
iteration closer to gradient descent.
The $\texttt{curve\_fit}$ routine utilizes a 
procedure $\texttt{resid}$
which takes a diffusion coefficient $D$ as input,
runs the Crank-Nicolson FD scheme \cref{ds:HE_CN} with this $D$,
and 
returns the residual vector
$\{U_\mathbf{j}^n - u_\mathbf{j}^n(D)\}$.
Note that even though the 
explicit forward Euler FD scheme \cref{ds:HE_FD} is faster, 
the implicit Crank-Nicolson scheme~\cref{ds:HE_CN} 
was chosen to avoid any stability issues 
that might arise in exploring the 
diffusion coefficient parameter space.

In general,
for a single-parameter least-squares optimization problem with a global minimum as in the 
present case, 
the LM algorithm is known to converge 
to a local minimum,
assuming the initial guess is sufficiently
close to the minimum~\cite{LM}.
In the present case,
we tested various initial guesses 
and convergence was obtained
in all cases, probably owing to the smooth
convex nature of the cost function shown in~\cref{fg:cost}.

%%%%%%%%%%%%%%%%%%%%%%%

\section{Diffusion coefficient estimates} \label{sec:est_res}

The procedure outlined in \Cref{sec:estim} 
was used to compute estimates 
$D_{\rm opt}(N)$
for the diffusion coefficient of 
argon in helium for several  
binning parameter values~$N = 20, 50, 100$.
The MD simulations were done with 
three different initial random seeds
and
the estimates of $D_{\rm opt}(N)$
were averaged.
The FD calculations used 
grid spacing $h = 1/N$.
\Cref{tb:diff_data} shows the results, where
column 1 gives the MD binning parameter~$N$,
column 2 gives the estimated $D_{\rm opt}(N)$ in cm$^2$/s,
column 3 gives the value of the cost function at $D_{\rm opt}(N)$ as defined in \cref{eq:cost},
column 4 gives the 95\% confidence interval for $D_{\rm opt}(N)$ as reported by $\texttt{curve\_fit}$.
The first three rows give the results 
of the estimation procedure 
for the different MD binning parameters,
and
the last row gives the 
experimentally measured diffusion coefficient
which was obtained using 
gas chromatography~\cite{Wasik1969-my}.

\begin{table}[htb] 
    \centering
    \begin{tabular}{|c|c|c|c|}
        \hline
MD binning parameter & $D_{\text{opt}}(N)$ [cm$^2$/s]  & $\texttt{cost}(D_{\rm opt}(N))$ & $95\%$ CI\\
        \hline
        $N=20$     & 0.9121        & 
         0.009870  & 1.48e-3 
        \\ \hline
        $N=50$     & 0.7948       & 0.008930 & 0.54e-3   
        \\ \hline
        $N=100$    & 0.8417       & 0.024214     & 0.46e-3    
        \\ \hline
       experiment~\cite{Wasik1969-my} & 0.7335 & & 8.70e-3 
        \\ \hline
    \end{tabular}
    \vskip 2.5pt
\caption{Results of 
diffusion coefficient estimation for
argon in helium,
column~1: MD binning parameter $N$,
column~2: optimal diffusion coefficient 
$D_{\rm opt}(N)$
averaged over three MD runs with different initial random seeds,
column~3: cost function 
$\texttt{cost}(D_{\rm opt}(N))$ 
as defined in~\cref{eq:cost},
column~4: 95\% confidence interval,
last row gives experimental value~$D_{\rm exp}$
of the diffusion coefficient~\cite{Wasik1969-my}. }
\label{tb:diff_data}
\end{table}

The following observations were made.
\begin{itemize}
\item 
Column~2 shows that
the estimated diffusion 
coefficient~$D_{\rm opt}(N)$ agrees best
with experiment for the intermediate
MD binning parameter $N=50$.
This can be understood as follows.
The smaller MD binning parameter $N=20$
makes the FD grid spacing $h=1/N$ larger, 
and 
then the FD result $u_{\bf j}^n$ 
is less accurate.
On the other hand,
the larger MD binning parameter $N=100$
makes the bins smaller
and hence they contain fewer atoms,
and
then the 
MD binning concentration $U_{\bf j}^n$
has more statistical noise.

\item 
Column 3 shows that the cost function
achieves its smallest value for $N=50$,
which is consistent with the previous point
that the best agreement between
model and experiment occurs for $N=50$.

\item
Column 4 shows that the model
confidence interval becomes narrower
with higher resolution;
however, the experimental confidence interval
is much wider,
which perhaps can be attributed to 
experimental measurement uncertainty.
\end{itemize}

\vspace*{5pt}
As shown in \cref{tb:diff_data}, 
although the estimated diffusion coefficient~$D_{\rm opt}(50)$ is close to the 
experimental value,
there is still a discrepancy
which may be attributed to 
several factors.
One factor is the choice of 
numerical parameters including the
number of atoms in the MD simulation~$N_p$,
binning parameter~$N$,
and
time steps~$k_{\rm MD},k_{\rm FD}$
in the MD and FD calculations.
A second factor is that the 
MD and FD calculations were initialized
differently and this could affect
the quality of the estimation results.
A third factor is that
the Lennard-Jones potential 
utilized in the MD simulations
is only an approximation to the
true interparticle interaction between the gas molecules.
A fourth factor is that the
calculations were carried out in 
two-dimensional space,
while there may be three-dimensional
effects in the experiment.
A final factor is the
statistical noise in the binned MD data 
which could be reduced with a 
convolution filter~\cite{larson_1997_hydrodynamics}.
These are all possible directions for future
investigation.

Finally,
we emphasize our choice to use a small sample of only three MD simulations at each $N$ value,
considering the relatively high MD runtime (up to 80 minutes).
We anticipate the utility of our method to be the extraction of an accurate diffusion coefficient estimate from a relatively small amount of MD data,
which may be applied to a efficiently simulate a larger continuum system.
Using more simulations may reduce the impact of the statistical noise observed in the binned MD concentration of \cref{fg:diff_results}.
However, 
we believe the reasonable agreement of the estimated and experimental diffusion coefficients justifies the method. 

%%%%%%%%%%%%%%%%%%

\section{Summary} \label{sec:sum}

This work studied the diffusion 
of a patch of argon gas in 
a background of helium gas
with the goal of determining the
diffusion coefficient of the argon atoms.
Two models were utilized,
(1) finite difference (FD) simulations of the 
continuum diffusion equation~\cref{eq:HE}
with a given diffusion coefficient $D$,
and
(2) molecular dynamics (MD) simulations
where the gas atoms were represented by
point particles interacting through
the Lennard-Jones potential~\cref{eq:6-12_pot}.
The MD simulations were taken as
the true physical representation of the
diffusion process.
To estimate the argon diffusion coefficient, 
the MD trajectories of the argon particles
were converted into
a binned concentration on the FD grid,
and
the optimal diffusion coefficient~$D_{\rm opt}$ was obtained
by minimizing the 
least squares difference between
the FD concentration
and
the binned MD concentration.
Comparison of the numerical results
with the experimental diffusion coefficient
suggests that our optimization
estimation approach was successful,
although there is still a discrepancy 
and
several directions for future investigation
were suggested.

The present work explored a system where an argon 
patch diffused into a background of helium.
Another interesting physical system is a mixture of polymers, 
where both phase-separation and diffusion dynamics are present.
Such systems may be modeled by the Allen-Cahn and Cahn-Hilliard equations,
\begin{equation}
\frac{\partial u}{\partial t}=M[\gamma\nabla^2u-u^3+u],
\quad
\frac{\partial u}{\partial t}=M\nabla^2[u^3-u-\gamma\nabla^2u],
\end{equation}
where 
$u$ is the order parameter
giving the purity of phases, 
$M$ is the mobility coefficient, 
and $\sqrt{\gamma}$ sets the 
width of the transition zone between 
the phases at equilibrium.
The parameter estimation algorithm developed in this article
may also be applied to find optimal values of $M$ and $\gamma$.
Some preliminary results for the phase separation equations
are contained in the REU project report~\cite{report}.

%%%%%%%%%%%%%%%%%%%%%%%%%%%%
\section*{Acknowledgments}

The author thanks the 
University of Michigan Department of Mathematics for hosting this project which was supported by National Science Foundation
grant DMS-2110767.

\bibliographystyle{siamplain}
\bibliography{bibliography.bib}

\end{document}